\input amstex\documentstyle{amsppt}  
\pagewidth{12.5cm}\pageheight{19cm}\magnification\magstep1
\topmatter
\title{A parametrization of unipotent representations}\endtitle
\author G. Lusztig\endauthor
\address{Department of Mathematics, M.I.T., Cambridge, MA 02139}\endaddress
\thanks{Supported by NSF grant DMS-1855773 and by a Simons Fellowship}\endthanks
\endtopmatter   
\document
\define\ovm{\overset\smile\to}

\define\Irr{\text{\rm Irr}}

\define\mpb{\medpagebreak}

\define\frl{\forall}
\define\pe{\perp}
\define\si{\sim}

\define\sqc{\sqcup}

\define\qua{\quad}

\define\tcl{\ti\cl}

\define\op{\oplus}
   
\define\part{\partial}
\define\emp{\emptyset}

\define\ra{\rangle}
\define\n{\notin}

\define\m{\mapsto}
\define\do{\dots}
\define\la{\langle}

\define\lra{\leftrightarrow}

\define\sub{\subset}    
\define\bxt{\boxtimes}
\define\T{\times}
\define\ti{\tilde}
\define\nl{\newline}
\redefine\i{^{-1}}

\define\un{\underline}
\define\ov{\overline}
\define\ot{\otimes}

\define\Hom{\text{\rm Hom}}

\define\Ind{\text{\rm Ind}}

\define\a{\alpha}
\redefine\b{\beta}

\define\g{\gamma}
\redefine\d{\delta}
\define\e{\epsilon}
\define\et{\eta}
\define\io{\iota}

\define\p{\pi}

\define\r{\rho}
\define\s{\sigma}
\redefine\t{\tau}
\define\th{\theta}
\define\k{\kappa}
\redefine\l{\lambda}
\define\z{\zeta}
\define\x{\xi}

\redefine\G{\Gamma}
\redefine\D{\Delta}

\define\Th{\Theta}
\redefine\L{\Lambda}

\redefine\aa{\bold a}

\define\hh{\bold h}

\redefine\ss{\bold s}

\define\BB{\bold B}
\define\CC{\bold C}

\define\EE{\bold E}
\define\FF{\bold F}

\define\HH{\bold H}

\define\NN{\bold N}

\define\PP{\bold P}

\define\RR{\bold R}
\define\SS{\bold S}

\define\ZZ{\bold Z}

\define\cc{\Cal C}

\define\cf{\Cal F}
\define\cg{\Cal G}
\define\ch{\Cal H}
\define\ci{\Cal I}

\define\cl{\Cal L}
\define\cm{\Cal M}

\define\cu{\Cal U}

\define\cx{\Cal X}

\define\fT{\frak T}

\define\tu{\ti u}

\define\tB{\ti B}

\define\tI{\ti I}

\define\sha{\sharp}

\define\bP{\bar P}

\define\tcf{\ti{\cf}}

\redefine\gg{\bold g}
\head Introduction\endhead
\subhead 0.1\endsubhead
Let $G$ be a simple algebraic group defined and split over a finite field $F_q$. Let
$\cu$ be the set of isomorphism classes of irreducible unipotent representations (over $\CC$)
of the finite group $G(F_q)$. Let $W$ be the Weyl group of $G$ and let $\Irr(W)$ be the
set of isomorphism classes of irreducible representations (over $\CC$) of $W$.
In \cite{L79} a partition of $\Irr(W)$ into families is described and in \cite{L84} a
partition $\cu=\sqc_c\cu_c$ of $\cu$ (with $c$ running over the families of $\Irr(W)$) is
introduced. Moreover, in \cite{L84, \S4} to any family $c$ we have associated a finite group
$\cg_c$ and a bijection
$$\cu_c\lra M(\cg_c).\tag a$$
Here, for any finite group $\G$, $M(\G)$ is the set of $\G$-conjugacy classes of pairs
$(x,\r)$ were $x\in\G$ and $\r$ is an irreducible representation (over $\CC$) of the
centralizer $Z_\G(x)$ of $x$ in $\G$; let
$\CC[M(\G)]$ (resp. $\NN[M(\G)]$) be the vector space 
of formal $\CC$-linear combinations of elements in $M(\G)$ and let $A_\G:\CC[M(\G)]@>>>\CC[M(\G)]$ be the non-abelian
Fourier transform of \cite{L79} (a linear isomorphism with square $1$).
Let $\NN[M(\G)]$ (resp. $\RR_{\ge0}[M(\G)]$) be the set of vectors of $\CC[M(\G)]$
which are linear combinations with coefficients in $\NN$ (resp. $\RR_{\ge0}$)
of elements in the basis $M(\G)$ of $\CC[M(\G)]$.
For any family $c$ we write $A_c$ instead of $A_{\cg_c}:\CC[M(\cg_c)]@>>>\CC[M(\cg_c)]$.

In this paper we are interested in defining for any family $c$ a basis $\b_c$ of the $\CC$-vector space
$\CC[M(\cg_c)]$ which has the properties (b)-(e) below.

(b) There is a unique bijection $\io:M(\cg_c)@>\si>>\b_c$, $(x,\r)\m\th(\x,\r)$ such that any
$(x,\r)\in M(\cg_c)$ appears with nonzero coefficient in $\io(\x,\r)$; this coefficient is
actually $1$.

(c) There exists a partial order $\le$ on $M(\cg_c)$ such that for any
$(x,\r)\in M(\cg_c)$ we have $\io(x,\r)=(x,\r$ plus an $\NN$-linear combination of
elements $(x',\r')\in M(\cg_c)$ which are $<(x,\r)$. In particular we have
$\b_c\sub\NN[M(\cg_c)]$.

(d) We have $A_c(\b_c)\sub\RR_{\ge0}[M(\cg_c)]$.

(e) The matrix of $A_c$ with respect to the basis $\b_c$ is upper triangular
(for some partial order on $\b_c$).
\nl
Such a basis has been constructed in \cite{L20} (where it is denoted by $\ti\BB_c$ and is
called the {\it new basis})
extending the results of \cite{L19}; property (e) for $\ti\BB_c$ is verified in
\cite{L20a}. The basis $\b_c$ considered in this paper is a somewhat modified form of
$\ti\BB_c$. If $G$ is of type $A,B,C$ we have $\b_c=\ti\BB_c$. If $G$ is of type $D$,
\cite{L20} contains in addition to the definition  of $\ti\BB_c$ in \cite{L20,\S1}, a
variant of that definition given in \cite{L20,\S2}; in this paper, $\b_c$ will be the same as the
variant in \cite{L20,\S2} (this is better suited for an extension of our results to nonsplit
even orthogonal groups, as will be shown elsewhere).
If $G$ is of exceptional type and $|c|\n\{4,17\}$ we have $\b_c=\ti\BB_c$.
If $G$ is of exceptional type and $|c|=4$ (resp. $|c|=17$), two (resp. three) elements in $\b_c$
differ from the corresponding elements in $\ti\BB_c$ (see the definition of $Prim(S_3), Prim(S_5)$
in 5.7); the reason for the change (at least in the case where $|c|=17$) is to make formulas
more symmetric.

While the definition of $\ti\BB_c$ in \cite{L20} for classical types was quite different from that for
exceptional types, in this paper the definition of $\b_c$ for classical types is similar to that for
exceptional types (the fact that such an approach is possible was stated without proof in \cite{L20}).
The basis $\b_c$ will be called the {\it second basis}
of $\CC[M(\cg_c)]$.

\subhead 0.2\endsubhead
Let $\PP$ be the set of pairs $(P,R)$ where $P$ is a parabolic 
subgroup of $G$ defined over $F_q$ with reductive quotient $\bP$ and $R$ is an irreducible
unipotent cuspidal representation over $\CC$ (up to isomorphism) of $\bP(F_q)$.
Let $\un\PP$ be the set of orbits of the conjugation action of $G(F_q)$ on $\PP$.
If $(P,R)\in\PP$, we can view $R$ as a representation of $P(F_q)$ and induce this from
$P(F_q)$ to $G(F_q)$. The irreducible representations of $G(F_q)$ appearing in this induced
representation form a subset $\cu^{P,R}$ of $\cu$ which depends only on the image
of $(P,R)$ in $\un\PP$. Hence for any $\gg\in\un\PP$ the subset $\cu^\gg$ of $\cu$
is well defined. We have $\cu=\sqc_{\gg\in\un\PP}\cu^\gg$ (see \cite{L78, 3.25}).

In this paper we fix a family $c$ of $\Irr(W)$.

We have a partition
$$\cu_c=\sqc_{\gg\in\un\PP}\cu_c^\gg\tag a$$
where $\cu_c^\gg=\cu_c\cap\cu^\gg$. Under the bijection 0.1(a), this corresponds to a partition
$$M(\cg_c)=\sqc_{\gg\in\un\PP}M(\cg_c)^\gg\tag b.$$

\subhead 0.3\endsubhead
In this paper we define a second partition of $\cu_c$ which in some sense is transversal to the
partition 0.2(a) of $\cu_c$ (see 0.6(b)).  
More precisely, we will describe a collection $\HH_c$ of subgroups of $\cg_c$ and a certain
subset $\ovm{\HH}_c$ of $\HH_c\T\HH_c$ such that for any $(H,H')\in\ovm{\HH}_c$, $H$ is a
normal subgroup of $H'$; we will also describe for any $(H,H')\in\ovm{\HH}_c$ a nonempty
subset $Prim(H,H')$ of $M(H'/H)$ such that

(a) $\Xi(c):=\{(H,H',e);(H,H')\in\ovm{\HH}_c,e\in Prim(H,H')\}$
\nl
is in canonical bijection $\Th:\Xi(c)@>\si>>M(\cg_c)$ 
(see 0.6(a)) with $M(\cg_c)$ (hence also with $\cu_c$). If we identify $M(\cg_c)$
with $\Xi(c)$ using this bijection, we obtain a surjective map $\a_c:M(\cg_c)@>>>\ovm{\HH}_c$
given by $(H,H',e)\m(H,H')$.
For $\hh\in\ovm{\HH}_c$ we set ${}^\hh M(\cg_c)=\a_c\i(\hh)$.
Under our identification $M(\cg_c)=\cu_c$, ${}^\hh M(\cg_c)$ becomes a subset
${}^\hh\cu_c$ of $\cu_c$; our second partition of $\cu_c$ is
$$\cu_c=\sqc_{\hh\in\ovm{\HH}_c}{}^\hh\cu_c.\tag c$$

\subhead 0.4\endsubhead
The set $\HH_c$ is in bijection with a subset $Con^+_c$ of $\NN[M(\cg_c)]$ defined as follows.

Let $\NN[c]$ be the set of formal $\NN$-linear combinations of elements in $c$.
In \cite{L84, \S4} an imbedding $c\sub M(\cg_c)$ is described; this induces an
imbedding $\NN[c]\sub\NN[M(\cg_c)]$. Let $Con_c$ be the set of constructible
representations of $W$ associated to $c$ in \cite{L82}; we view the elements of $Con_c$ as
elements of $\NN[c]$ hence, using the imbedding above, as elements of $\NN[M(\cg_c)]$. As
shown in \cite{L86},

(a) {\it the representations in $Con_c$ are precisely the representations of $W$
carried by the left cells of $W$ contained in the two-sided cell attached to $c$.}
\nl
We set
$$r_!=\sum_{(x,\r)\in M(\cg_c);x=1}\dim(\r)(x,\r)\in\NN[M(\cg_c)],$$
(see also 0.5(a)),
$$Con^+_c=Con_c\cup\{r_!\}\sub\NN[M(\cg_c)].$$
When $G$ is of classical type (but not in general) we have $r_!\in Con_c$ hence $Con^+_c=Con_c$.

\subhead 0.5\endsubhead
We explain how the subgroups in $\HH_c$ are attached to the various elements of $Con^+_c$.

Let $\G$ be a finite group and let $H$ be a subgroup of $\G$.
Following \cite{L84, p.312} we define a linear map $i_{H,\G}:\CC[M(H)]@>>>\CC[M(\G)]$ by 
$$(x,\s)\m\sum_\r(\r:\Ind_{Z_H(x)}^{Z_\G(x)}(\s))(x,\r).$$
where $\r$ runs over the irreducible representations of $Z_\G(x)$
up to isomorphism and $:$ denotes multiplicity. 

Assume now that $H$ is a normal subgroup of $\G$; let $\p:\G@>>>\G/H$ be the canonical map.
Following {\it loc.cit.} we define a linear map $\p_{H,\G}:\CC[M(\G/H)]@>>>\CC[M(\G)]$ by 
$$(x,\s)\m
\sum_{y\in\p\i(x)}\sum_{\t\in\Irr(Z_\G(y))}|Z_\G(y)||Z_{\G/H}(x)|\i|H|\i(\t:\s)(y,\t)$$    
where $\t$ runs over the irreducible representations of $Z_\G(y)$ up to isomorphism and
$(\t:\s)$ denotes the 
multiplicity of $\t$ in $\s$ viewed as a representation of $Z_\G(y)$ via the obvious homomorphism
$Z_\G(y)@>>>Z_{\G/H}(x)$. 
Now let $H\sub H'$ be two subgroups of $\G$ such that $H$ is normal in $H'$.
We define a linear map $\ss_{H,H';\G}:\CC[M(H'/H)]@>>>\CC[M(\G)]$ by $f\m i_{H',\G}(\p_{H,H'}(f))$. Note that

(a) $r_!=i_{\{1\},\G}(1,1)=\ss_{\{1\},\{1\};\G}(1,1)$.

In \cite{L87} to each $r\in Con_c$ we have attached a conjugacy class $[r]$ of subgroups of
$\cg_c$. Its main property is that for any $H\in[r]$ we have

(b) $r=\ss_{H,H;\cg_c}(1,1)$.
\nl
When $r=r_!$ we define $[r]$ to consist of $\{1\}$; then (b) continues to hold (see (a)).
Thus $[r]$ is defined for any $r\in Con^+_c$.
For each $r$ we choose a specific $H(r)\in[r]$ as follows.
When $G$ is of classical type, $\cg_c$ is abelian and $H(r)$ is the unique subgroup in $[r]$;
when $G$ is of exceptional type, the subgroups $H(r)$ are described in \S5.

By definition, we have
$$\HH_c=\{H(r);r\in Con^+_c\}.$$
Note that $r\m H(r)$ is a bijection $Con^+_c@>\si>>\HH_c$.

\subhead 0.6\endsubhead
The subset $\ovm{\HH}_c$ of $\HH_c\T\HH_c$ and the subsets $Prim(H,H')$ of $M(H'/H)$ (for
$(H,H')\in\ovm{\HH}_c$) in 0.3 are described:

(i) in 3.7 when $G$ is isogenous to a symplectic or odd special orthogonal group;

(ii) in 3.8 when $G$ is isogenous to an even special orthogonal group;

(iii) in \S5 when $G$ is of exceptional type.
\nl
(When $G$ is of type $A$, $\HH_c$ consists of $\{1\}$, 
 $\ovm{\HH}_c$ consists of $(\{1\},\{1\})$ and $\Xi_c$ consists of one element.)

\mpb

We now state one of our main results.

(a) {\it There is a unique bijection $\Xi_c@>\Th>>M(\cg_c)$ (notation of 0.3(a)) such that for
any $(H,H',e)\in\Xi_c$, the element $\Th(H,H',e)\in M(\cg_c)$ appears with nonzero coefficient
in $\ss_{H,H';\cg_c}(e)\in\CC[M(\cg_c)]$.}
\nl
In view of 0.1(a), this can be viewed as a parametrization of $\cu_c$.

The proof of (a) when $G$ is as in (i) is given in 3.7; it is based on
the results on the new basis in \cite{L20}.
The proof of (a) when $G$ is as in (ii) is given in 3.8; it is based on
an extension of the results on the new basis in \cite{L20}.
The proof of (a) when $G$ is as in (iii) is given in \S5; a result close to (a)
was already known in this case from \cite{L20}.
When $G$ is of type $A$, the statement (a) is obvious.

Here is one of the main properties of the bijection $\Th$.

(b) {\it The partition 0.3(b) of $\cu_c$ is transversal to the partition 0.2(a) of $\cu_c$, in
the sense that for any $\gg\in\un\PP$ and any $\hh\in\ovm{\HH}_c$ we have
$|{}^\hh\cu_c\cap\cu_c^\gg|\le1$.}

From the proof of (a) and the results of \cite{L20} we see that
$\b_c:=\{\ss_{H,H';\cg_c}(e);(H,H',e)\in\Xi_c\}$ (see 0.6(a)) is a basis of $\CC[M(\cg_c)]$
(which could be called the {\it second basis} of $\CC[M(\cg_c)]$) in bijection with $\Xi_c$ and which
satisfies 0.1(b)-(d); it is also in canonical bijection

$\ss_{H,H';\cg_c}(e)\m\Th(H,H',e)$
\nl
with the obvious basis of $\CC[M(\cg_c)]$) consisting of
the  various elements of $M(\cg_c)$. From \cite{L20a} or some variant of it (in type $D$)
we see that $\b_c$ satisfies 0.1(e).

For any $(x,\r)\in M(\cg_c)^\gg$ (see 0.2(b))
let $g_{x,\r}$ be the element of the second basis corresponding as
above to $(x,\r)$; we define an element $f_{x,\r}\in\CC[M(\cg_c)]$ by the following
requirement: for any $(x',\r')\in M(\cg_c)^{\gg'}$ the coefficient of
$(x',\r')$ in $f_{x,\r}$ is equal
to the coefficient of $(x',\r')$ in $g_{x,\r}$ if $\gg=\gg'$ and is equal to $0$ if $\gg\ne\gg'$.

We have the following result.

(c) {\it $\{f_{x,\r};(x,\r)\in M(\cg_c)\}$ is a basis of $\CC[M(\cg_c)]$.}
\nl
(This could be called the {\it third} basis of $\CC[M(\cg_c)]$.)

The proof in the case 0.6(i) is given in \S7. The proof in the case 0.6(ii) is similar.
The proof in the case 0.6(iii) is obtained by examining the tables.

\subhead 0.7\endsubhead
Let $\gg\in\un\PP$. Let
$$\ovm{\HH}_c^\gg=\{\hh\in\ovm{\HH}_c;|{}^\hh\cu_c\cap\cu_c^\gg|=1\}$$.

From 0.6(a) we see that we have a bijection

(a) $\ovm{\HH}_c^\gg@>>>\cu_c^\gg$
\nl
which to any $\hh\in\ovm{\HH}_c^\gg$ associates the unique element in ${}^\hh\cu_c\cap\cu_c^\gg$.
We thus obtain the following statement which is one of the main results of this paper.

(b) {\it $\cu_c^\gg$ is in natural bijection with a subset of $\ovm{\HH}_c$, hence with a
subset of $\HH_c\T\HH_c$.}

\subhead 0.8\endsubhead
A key role in this paper is played by the subset $\ovm{\HH}_c$ of $\HH_c\T\HH_c$.
As mentioned above, $\HH_c$ can be identified with $Con^+_c$. Therefore 
$\ovm{\HH}_c$ can be identified with a subset of $Con^+_c\T Con^+_c$ which we can denote by
$\ovm{Con^+_c}$ and also with the image of this set in
$(Con^+_c\T Con^+_c)_{unord}$ (unordered pairs in $Con^+_c\T Con^+_c$).

We will show elsewhere that $\ovm{Con^+_c}$ admits a direct (inductive) definition.

We can restate 0.7(b) to say that, if $\gg\in\un\PP$, then

(a) {\it $\cu_c^\gg$ is in natural bijection with a subset of $\ovm{Con^+_c}$, hence with a
subset of $(Con^+_c\T Con^+_c)_{unord}$.}

\subhead 0.9\endsubhead
{\it Erratum to \cite{L20}.} In the first displayed equality of 3.2 replace $Z_\G(z)$ by $Z_\G(x)$.
In line 8 of 3.7 replace $H_{221}$ by $H_{211}$.

\subhead 0.10\endsubhead
{\it Notation.} $\FF$ denotes the field with two elements. For $a,b\in\ZZ$ we set
$[a,b]=\{z\in\ZZ;a\le z\le b\}$. We write $a\ll b$ instead of $b-a\ge2$. 

\head Contents\endhead
1. The set $\pmb\cf(V_D)$.

2. The bijection $\mu:\un{\pmb\cf}(V_D)@>\si>>\pmb\cf(V_D)$.

3. From  $\SS_D$ to $V_D$.

4. Tables in classical types.

5. Exceptional types.

6. The partition $\cu_c=\sqc_{\gg\in\un\PP}\cu_c^\gg$.

7. The third basis of $\CC[M(\cg_c)]$.

\head 1. The set $\pmb\cf(V_D)$\endhead
\subhead 1.1\endsubhead
We fix $D\in\NN$. Let $V_D$ be the $\FF$-vector space with basis $\{e_i;i\in[1,D]\}$. 
When $D'\in[0,D]$ we identify $V_{D'}$ with the subspace of $V_D$ with basis $\{e_i,i\in[1,D']\}$.
When $D\ge2$, for any $i\in[1,D]$ there is a unique linear map $T_i:V_{D-2}@>>>V_D$ such that

$T_i(e_k)=e_k$ if $k\in[1,i-2]$,

$T_i(e_k)=e_{k+2}$ if $k\in[i,D-2]$,

$T_i(e_{i-1})=e_{i-1}+e_i+e_{i+1}$ if $1<i<D$.
\nl
This map is injective.

If $v_1,v_2,\do,v_k$ are vectors in $V_D$ we denote by $\la v_1,v_2,\do,v_k\ra$ the
subspace of $V_D$ generated by $v_1,v_2,\do,v_k$.

\subhead 1.2\endsubhead
A subset $I\sub\ZZ$ is said to be an {\it interval} if it is of the form
$[a,b]$ for some $a\le b$ in $\ZZ$. Let $\ci_D$ be the set of intervals contained in $[1,D]$.
We say that $I=[a,b],I'=[a',b']$ in $\ci_D$ are non-touching if $b\ll a'$ or if $b'\ll a$.
For $I=[a,b],I'=[a',b']$ in $\ci_D$ we write $I\prec I'$ if $a'<a\le b<b'$.

Let $R_D$ be the set whose elements are the subsets of $\ci_D$.
We define a subset $\SS^{prim}_D$ of $R_D$ as follows.
When $D$ is even, $\SS^{prim}_D$ consists of $\emp$ and of
$$\{[1,D],[2,D-1],\do,[k,D+1-k]\}\tag a$$

for various $k\in[1,D/2]$.
When $D$ is odd, $\SS^{prim}_D$ consists of $\emp$, of
$$\{[1,D-1],[2,D-2],\do,[k,D-k]\}\tag b$$

and of
$$\{[2,D],[3,D-1],\do,[k+1,D+1-k]\}\tag c$$

for various odd $k\in[1,(D-1)/2]$. 

For example, if $D=2$, $\SS^{prim}_D$ consists of $\emp,\{[1,2]\}$.

If $D=4$, $\SS^{prim}_D$ consists of $\emp,\{[1,4]\},\{[1,4],[2,3]\}$.

If $D=6$, $\SS^{prim}_D$ consists of $\emp,\{[1,6]\},\{[1,6],[2,5]\}$,
$\{[1,6],[2,5],[3,4]\}$.   

If $D=3$, $\SS^{prim}_D$ consists of $\emp,\{[1,2]\},\{[2,3]\}$.

If $D=5$, $\SS^{prim}_D$ consists of $\emp,\{[1,4]\},\{[2,5]\}$.

If $D=7$, $\SS^{prim}_D$ consists of

$\emp,\{[1,6]\},\{[2,7]\}$, $\{[1,6],[2,5],[3,4]\},\{[2,7],[3,6],[4,5]\}$.

The elements of $\SS_D^{prim}$ are said to be {\it primitive.}

When $D\ge2$ and $i\in[1,D]$ we define an (injective) map $\x_i:\ci_{D-2}@>>>\ci_D$ by

$\x_i([a',b'])=[a'+2,b'+2]$ if $i\le a'$,

$\x_i([a',b'])=[a',b']$ if $i\ge b'+2$,

$\x_i([a',b'])=[a',b'+2]$ if $a'<i<b'+2$.
\nl
We define $t_i:R_{D-2}@>>>R_D$ by $B'\m\{\x_i(I');I'\in B'\}\sqc\{i\}$ (see \cite{L20, 1.1}).

We define a subset $\SS_D$ of $R_D$ by induction on $D$ as follows. 

If $D=0$, $\SS_D$ consists of $\emp\in R_D$. If $D=1$, $\SS_D$ consists of $\emp$ and of
$\{1\}\in R_D$. If $D\ge2$, a subset $B$ of $\ci_D$ is in $\SS_D$ if either
$B\in\SS^{prim}_D$ or if there exists $i\in[1,D]$ and $B'\in\SS_{D-2}$ such that $B=t_i(B')$.
(When $D$ is even this definition appears in \cite{L20}.)

We note the following result (when $D$ is even this appears in \cite{L20}).

(d) {\it If $B\in\SS_D$ and $I\in B,\tI\in B$, then either $I=\tI$ or
$I,I'$ are non-touching or $I\prec\tI$ or $\tI\prec I$.}
\nl
We use induction on $D$.
If $B$ is primitive then (d) is obvious. We now assume that $B$ is not primitive.
Then $D\ge2$ and there exists $i\in[1,D]$ and $B'\in\SS_{D-2}$ such that $B=t_i(B')$.
By the induction hypothesis, (d) holds when $B$ is replaced by $B'$.
It follows immediately that (d) holds for $B$.

\mpb

For any $I\sub[1,D]$ we set $e_I=\sum_{i\in I}e_i\in V_D$.
For $B\in\SS_D$ let $\EE_B$ be the subspace of $V_D$ spanned by $\{e_I;I\in B\}$.
Let $\pmb\cf(V_D)$ be the set of subspaces of $V_D$ of the form $\EE_B$ for some $B\in\SS_D$.
We have the following result (when $D$ is even this appears in \cite{L20}).

(e) {\it If $B\in\SS_D$ then $\{e_I;I\in B\}$ is a basis of $\EE_B$.} 
\nl
We can assume that $D\ge2$. Assume that there exists a nonempty subset $\cx\sub B$ such that
$\sum_{I\in\cx}e_I=0$. Let $I_1=[a,b]\in\cx$ be such that $|I_1|$ is maximum.
If $I_2\in\cx$, $I_2\ne I_1$ and $a\in I_2$, then $a\in I_1\cap I_2$ so that by (d) we have
$I_2\prec I_1$ or $I_1\prec I_2$; now $I_2\prec I_1$ contradicts $a\in I_1,a\in I_2$
and $I_1\prec I_2$ contradicts the maximality of $|I_1|$. We see that 
if $I_2\in\cx$, $I_2\ne I_1$, then $a\n I_2$ (and similarlly $b\n I_2$).
It follows that the coefficient of $e_a$ (and that of $e_b$) in 
$\sum_{I\in\cx}e_I$ is $1$. But the last sum is zero. This contradiction proves (e).

We have the following result.

(f) {\it If $B\in\SS_D$ and $J\in\ci_D$ is such that $e_J\in\EE_B$ then $J\in B$.}
\nl
Assume that $J\n B$. We can find a nonempty subset $\cx\sub B$ such that
$\sum_{I\in\cx}e_I=e_J$. Let $I_1=[a,b]\in\cx$ be such that $|I_1|$ is maximum.
As in the proof of (e) we see that the coefficient of $e_a$ (and that of $e_b$) in 
$\sum_{I\in\cx}e_I$ is $1$. Since the last sum is equal to $e_J$ it follows that $a\in J,b\in J$.
Since $J$ is an interval we have $[a,b]\sub J$, Since $J$ is an interval different from $[a,b]$
we must have $a-1\in J$ or $b+1\in J$.
Assume first that $a-1\in J$. Since $\sum_{I\in\cx}e_I=e_J$ we can find
$I_2\in\cx$ such that $a-1\in I_2$.
We have $I_2\not\sub I_1$ (since $a-1\in I_2,a\n I_1$); we cannot have $I_1\prec I_2$
(this would contradict the
maximality of $|I_1|$). Using (d) we deduce that $I_1,I_2$ are non-touching, but
this contradicts $a-1\in I_2,I_1=[a,b]$. We see that $a-1\in J$ leads to a contradiction.
Similarly $b+1\in J$ leads to a contradiction.
This contradiction proves (f).

(g) {\it The map $\SS_D@>>>\pmb\cf(V_D)$, $B\m\EE_B$ is a bijection.}
\nl
It is enough to show that this map is injective. Let $B\in\SS_D,B'\in\SS_D$ be such that
$\EE_B=\EE_{B'}$. By (f) we have

$B=\{J\in\ci_D;e_J\in\EE_B\}=\{J\in\ci_D;e_J\in\EE_{B'}\}=B'$.
\nl
We see that $B=B'$.
This proves (g).

In the case where $D$ is even, (g) appears also in \cite{L20} but the present proof is
simpler.

\subhead 1.3\endsubhead
The set $\pmb\cf(V_D)$ has an alternative definition (by induction on $D$).
When $D\ge2$ we define a set $P(V_D)$  of subspaces of $V_D$ (said to be
{\it primitive} subspaces) as follows. If $D$ is even, $P(V_D)$ consists of
$\{0\}$ and of the subspaces $\EE_B$ with $B$ as in 1.2(a) and $k\in[1,D/2]$.
If $D$ is odd, $P(V_D)$ consists of $\{0\}$ and of the subspaces $\EE_B$ with $B$ as in
1.2(b) or 1.2(c) and odd $k\in[1,(D-1)/2]$.
If $D=0$, $\pmb\cf(V_D)$ consists of the subspace $\{0\}$.
If $D=1$, $\pmb\cf(V_D)$ consists of the subspace $\{0\}$ and of $V$.
If $D\ge2$, a subspace $\EE$ of $V_D$ is in $\pmb\cf(V_D)$ if it is either
primitive or if there exists $i\in[1,D]$ and $\EE'\in\pmb\cf(V_{D-2})$ such that $\EE=T_i(\EE')\op\FF e_i$. 

\subhead 1.4\endsubhead
For $\d\in\{0,1\}$ we set $\ci_D^\d=\{I\in\ci_D;|I|=\d\mod2\}$.
For $B\in\SS_D$ we set $B^\d=B\cap\ci_D^\d$. Let
$$S_D=\{B\in\SS_D;B=B^1\}.$$
We can also define $S_D$ by induction on $D$ as follows.
We have $S_0=\{\emp\}$; $S_1$ consists of $\{\emp\}$ and $\{1\}$. If $D\ge2$ then 
a subset $B$ of $\ci_D$ is in $S_D$ if either $B=\emp$ or if
there exists $i\in[1,D]$ and $B'\in S_{D-2}$ such that $B=t_i(B')$.

We define a collection $\cf(V_D)$ of subspaces of $V_D$ by induction on $D$ as follows. 
If $D=0$, $\cf(V_D)$ consists of $\{0\}$. If $D=1$, $\cf(V_D)$ consists of $\{0\}$ and $V_D$.
If $D\ge2$, a subspace $E$ of $V_D$ is said to be in $\cf(V_D)$ if either $E=0$ or if there
exists $i\in[1,D]$ and $E'\in\cf(V_{D-2})$ such
that $E=T_i(E')\op\FF e_i$. We have $\cf(V_D)\sub\pmb\cf(V_D)$.
Note that, for $\EE\in\pmb\cf(V_D)$ with $\EE=\EE_B$, $B\in\SS_D$ we have
$\EE\in\cf(V_D)$ if and only $B\in S_D$.

\subhead 1.5\endsubhead
Let $\d\in\{0,1\}$. Let $\ZZ^\d=\d+2\ZZ$. For any $I\in\ci_D$ we set $I^\d=I\cap\ZZ^\d$;
we have $I=I^0\sqc I^1$.
Let $V_D^\d=\la e_i;i\in[1,D]^\d\ra$. We have $V_D=V_D^0\op V_D^1$.

Assume that $D\ge2$, $i\in[1,D]$. There is a unique linear map
$T_i^\d:V_{D-2}^\d@>>>V_D^\d$ such that

$T_i^\d(e_k)=e_k$ if $k\in[1,i-2]^\d$;

$T_i^\d(e_k)=e_{k+2}$  if $k\in[i,D-2]^\d$;

$T_i^\d(e_{i-1})=e_{i-1}+e_{i+1}$ if $i\in[2,D-1]^{1-\d}$.
\nl
Note that for $x\in V_{D-2}^0,y\in V_{D-2}^1$ we have

(a) $T_i(x+y)=T_i^0(x)+T_i^1(y)\mod\FF e_i$.

\subhead 1.6\endsubhead
We define a collection $\cc(V_D^\d)$ of subspaces of $V_D^\d$ by induction on $D$ as
follows. If $D=0$, $\cc(V_D^\d)$ consists of $\{0\}$. If $D=1$, $\d=0$, $\cc(V_D^\d)$ consists
of $\{0\}$. If $D=1$, $\d=1$, $\cc(V_D^\d)$ consists of $\{0\}$ and $V_D^\d$.
If $D\ge2$, a subspace $\cl$ of $V_D^\d$ is said to be in $\cc(V_D^\d)$ if either 

there exists $i\in[1,D]^\d$ and $\cl'\in\cc(V_{D-2}^\d)$ such that $\cl=T_i^\d(\cl')\op\FF e_i$, or

there exists $i\in[1,D]^{1-\d}$ and $\cl'\in\cc(V_{D-2}^\d)$ such that $\cl=T_i^\d(\cl')$.
\nl
For $E\in\cf(V_D)$ we set $E^\d=E\cap V_D^\d$. Then

(a) $E^\d\in\cc(V_D^\d)$, $E=E^0\op E^1$.
\nl
When $D$ is even this is shown in \cite{L19, 2.2(c), 2.3(b)}. The case where $D$ is odd
is similar.

\subhead 1.7\endsubhead
We define by an induction on $D$ a collection $\tcf(V_D)$ of pairs $(\cm^0,\cm^1)$ where
$\cm^0,\cm^1$ are subspaces of $V_D^0,V_D^1$ respectively.
If $D=0$, $\tcf(V_D)$ consists of $(\{0\},\{0\})$.
If $D=1$, $\tcf(V_D)$ consists of $(\{0\},\{0\})$ and $(\{0\},V_D^1)$.
If $D\ge2$ a pair $(\cm^0,\cm^1)$ of subspaces of $V_D^0,V_D^1$ is said to be in $\tcf(V_D)$ if 
either $(\cm^0,\cm^1)=(\{0\},\{0\})$ or there exists $i\in[1,D]$ and
$(\cm'{}^0,\cm'{}^1)\in\tcf(V_{D-2})$ such that

if $\d=i\mod2$ then $\cm^\d=T_i^\d(\cm'{}^\d)\op\FF e_i$;

if $\d=i+1\mod2$, then $\cm^\d=T_i^\d(\cm'{}^\d)$.
\nl
Using 1.5(a) and the definitions we see that

(a) $E\m(E^0,E^1)$ is a bijection $\cf(V_D)@>\si>>\tcf(V_D)$. Moreover, if
$(\cm^0,\cm^1)\in\tcf(V_D)$, then $\cm^\d\in\cc(V_D^\d)$ for $\d=0,1$.

\subhead 1.8\endsubhead
There is a unique symplectic form $(,):V_D\T V_D@>>>\FF$ such that for $i,j$ in $[1,D]$ we
have $(e_i,e_j)=1$ if $i-j=\pm1$, $(e_i,e_j)=0$ if $i-j\ne\pm1$. This form is nondegenerate if
$D$ is even and has a one dimensional radical spanned by
$$\et_D:=e_1+e_3+e_5+\do+e_D$$
if $D$ is odd. The next result follows from 1.2(d).

(a) {\it If $B\in\SS_D$, then $(,)$ is identically zero on $\EE_B\T\EE_B$,}
\nl
If $D\ge2$ then for $i\in[1,D]$ we have

(b) $(x,y)=(T_i(x),T_i(y))$ for any $x,y$ in $V_{D-2}$,

(c) $(T_i^0(x),T_i^1(y))=(x,y)$ for any $x\in V_{D-2}^0,y\in V_{D-2}^1$.
\nl
For any subspace $Z$ of $V_D$ we set $Z^\pe=\{x\in V_D;(x,y)=0\qua\frl y\in Z\}$.
When $Z\sub V_D^{1-\d}$, we set $Z^!=\{x\in V_D^\d;(x,y)=0\qua\frl y\in Z\}=Z^\pe\cap V_D^\d$.
Let

$\cf_*(V_D)=\{E\in\cf(V_D);\dim(E)=D/2\}$ (if $D$ is even),

$\cf_*(V_D)=\{E\in\cf(V_D);\dim(E)=(D+1)/2\}$ (if $D$ is odd).
\nl
We have the following result.

(d) {\it Assume that $D$ is even or that $D$ is odd and $\d=0$. If $\cl\in\cc(V_D^\d)$,
then $\cl^!\in\cc(V_D^{1-\d})$ and $\cl\op\cl^!\in\cf(V_D)$.
Moreover, $\cl\m\cl\op\cl^!$ is a bijection $\cc(V_D^\d)@>>>\cf_*(V_D)$.}
\nl
When $D$ is even this is proved in \cite{L19}. The proof for $D$ odd is similar.

We show:

(e) {\it If $D$ is odd, then $|\cf_*(V_D)|=|\cf_*(V_{D-1})|$.}
\nl
Using (d) we see that it is enough to show that $|\cc(V_D^0)|=|\cc(V_{D-1}^0)|$.
But from the definitions we actually have $V_D^0=V_{D-1}^0$ and $\cc(V_D^0)=\cc(V_{D-1}^0)$.

\subhead 1.9\endsubhead
Let $\d\in\{0,1\}$. Assume that $D$ is even or that $D$ is odd and $\d=1$.
We define a collection $\ovm{\cc}(V_D^\d)$ of pairs $(\cl\sub\tcl)$ of subspaces of $V_D^\d$.
Namely, $\ovm{\cc}(V_D^\d)$ consists of all pairs $(\cm^\d\sub(\cm^{1-\d})^!)$ with
$(\cm^0,\cm^1)\in\tcf(V_D)$. 
(We use that $\cm^0\op\cm^1$ is an isotropic subspace of $V_D$, see 1.7(a), 1.8(a).)

From 1.7(a), 1.8(d) we see that

(a) $\ovm{\cc}(V_D^\d)\sub\cc(V_D^\d)\T\cc(V_D^\d)$.
\nl
From 1.7(a) we see that we have a bijection

(b) $\cf(V_D)@>\si>>\ovm{\cc}(V_D^\d)$, $E\m(E^\d\sub(E^{1-\d})^!)$.
\nl
The inverse of the map (b) is $(\cl\sub\tcl)\m\cl\op\tcl^!$.

\subhead 1.10\endsubhead
In this subsection we assume that $D\ge1$. We show:

(a) {\it If $B\in S_{D-1}$ then $B\in S_D$.}
\nl
This makes sense since $\ci_{D-1}\sub\ci_D$. We argue by induction on $D$. If $D\le2$,
the result is obvious. Assume now that $D\ge3$. The maps $R_{D-3}@>>>R_{D-1}$
analogous to $t_i:R_{D-2}@>>>R_D$ (see 1.2) are denoted by $\ti t_i$;
they are defined for $i\in[1,D-1]$.
If $B=\emp$ the result is obvious. Thus we can assume that $B\ne\emp$. We can find
$B'\in S_{D-3}$ and $i\in[1,D-1]$ such that $B=\ti t_i(B')$.
From the definitions we have $\ti t_i(B')=t_i(B')$ so that $B=t_i(B')$.
By the induction hypothesis we have $B'\in S_{D-2}$ so that $B\in S_D$.
This proves (a).

\mpb

Note that $\SS_{D-1}$ is not in general contained in $\SS_D$.

\subhead 1.11\endsubhead
We now give an alternative (non-inductive) definition fo $S_D$.

(a) {\it Let $B\in R_D$ be such that $|I|$ is odd for any $I\in B$.
Then $B\in S_D$ if and only if $B$ satisfies properties $(P_0),(P_1)$ below.}

{\it $(P_0)$. If $I\in B,\tI\in B$, then either $I=\tI$, or
$I,I'$ are non-touching, or $I\prec\tI$, or $\tI\prec I$.}

{\it $(P_1)$. If $[a,b]\in B$ and $c\in\NN$ satisfy $a<c<b$, $c\ne a\mod2$, then there exists
$[a_1,b_1]\in B$ such that $a<a_1\le c\le b_1<b$.}
\nl
When $D$ is even this is proved in \cite{L19,1.3(c)}. The proof in the case where $D$ is odd is
similar. (The fact that any $B\in S_D$ satisfies $(P_0)$ is also contained in 1.2(d).)

Note that (a) provides an alternative proof of the inclusion $S_{D-1}\sub S_D$ in 1.10.

\head 2. The bijection $\mu:\un{\pmb\cf}(V_D)@>\si>>\pmb\cf(V_D)$\endhead
\subhead 2.1\endsubhead
Until the end of 2.6 we fix $B\in S_D$. Let $E=\EE_B\in\cf(V_D)$.
We have $|B|\le D/2$ if $D$ is even, $|B|\le(D+1)/2$ if $D$ is odd. We set $\s\in\cup_{I\in B}$.
We can write uniquely 
$$\s=[a_1,b_1]\cup[a_2,b_2]\cup\do\cup[a_s,b_s]$$
where $1\le a_1\le b_1\ll a_2\le b_2\ll\do\ll a_s\le b_s\le D$,
$[a_1,b_1]\in B,[a_2,b_2]\in B,\do,[a_s,b_s]\in B$. Note that
$c_i:=b_i-a_i\in2\NN$ for $i=1,\do,s$. We set $b_0=-1,a_{s+1}=D+2$. Note that
$$a_i-b_{i-1}\ge2\text{ for }i\in[1,s+1].$$
We have
$$|B|=\sum_{1\in[1,s]}(c_i+2)/2.\tag a$$
This can be proved by induction on $D$ in the same way as \cite{L20, 1.3(g)}.

In the remainder of this subsection we assume that $D$ is odd. We show:

(b) {\it If $\et_D\in E$ (notation of 1.8) then $|B|=(D+1)/2$.}
\nl
By our assumption, $\cup_{I\in B}$ contains $\{1,3,\do,D\}$. It follows that

$a_1,b_1,a_2,b_2,\do,a_s,b_s$
\nl
are all odd and
$$a_1=1,a_2=b_1+2,a_3=b_2+2,\do,a_s=b_{s-1}+2,b_s=D.$$
\nl
Using (a) we deduce
$$2|B|=\sum_{1\in[1,s]}(b_i-a_i+2)=-1-2-2-\do-2+2s+D=-1-2(s-1)+2s+D=D+1$$
and (b) is proved.

We show

(c) {\it If $|B|\le(D-1)/2$ then $\dim(E^\pe)=D-\dim(E)$.}
\nl
Since the radical of $(,)$ is spanned by $\et_D$, it is enough to show that
$\et_D\n E$. This follows from (b).

We show:

(d) {\it If $|B|=(D+1)/2$, then $\et_D\in E$.}
\nl
We argue by induction on $D$. If $D=1$ we have $E=V_1$ so that there is nothing to prove.
Assume now that $D\ge3$. We can find $i\in[1,D]$ and $E'\in\cf(V_{D-2})$ such that
$E=T_i(E')\op\FF e_i$. We have $\dim(E')=(D-1)/2$ hence by the induction hypothesis we have
$\et_{D-2}\in E'$. From the definitions we have $T_i(\et_{D-2})=\et_D\mod \FF e_i$ so that
$\et_D\in E$. This proves (d).

We show 

(e) {\it If $|B|=(D-1)/2$, then $\et_D\n E$ and $\FF\et_D\op E\in\cf(V_D)$.}
\nl
The first assertion follows from (b). For the second assertion we argue by induction on $D$.
If $D=1$, we have $E=0$ and there is nothing to prove.
Assume now that $D\ge3$. We can find $i\in[1,D]$ and $E'\in\cf(V_{D-2})$ such that
$E=T_i(E')\op\FF e_i$. We have $\dim(E')=(D-3)/2$ hence by the induction hypothesis we have
$\FF\et_{D-2}+E'\in\cf(V_{D-2})$. It follows that
$\FF T_i(\et_{D-2})+T_i(E')+\FF e_i\in\cf(V_D)$ hence (as in the proof of (d)) we have
$\FF\et_D+T_i(E')+\FF e_i\in\cf(V_D)$, that is $\FF\et_D+E\in\cf(V_D)$. This proves (e).

We show:

(f) {\it If $\dim(E)=(D+1)/2$, then there is a unique $E'\in\cf(V_{D-1})$ such that
$\dim(E')=(D-1)/2$ and $E=\FF\et_D+E'$.}
\nl
We define $B'$ to be $B$ from which the unique interval of $B$ containing $D$ is removed.
We have $B'\in S_{D-1}$ (viewed as a subset of $S_D$, see 1.10).
Let $E'\in\cf(V_{D-1})$ be the subspace of $V_{D-1}$ with basis $\{e_I;I\in B'\}$. We have
$\dim(E')=(D-1)/2$ and $E'\sub E$. By (d) we have $\et_D\in E$. It follows that
$\FF\et_D\op E'\sub E$ (the sum is direct by (e)). Since $\dim(\FF\et_D\op E')=\dim(E)=(D+1)/2$,
it follows that $\FF\et_D\op E'=E$. This proves the existence in (f). We now define a map
$$\{E'_1\in\cf(V_{D-1});\dim(E'_1)=(D-1)/2\}@>>>\{E_1\in\cf(V_D);\dim(E_1)=(D+1)/2\}$$
by $E'_1\m\FF\et_D+E'_1$. This map is well defined by (e) and is surjective by the first part
of the proof. Our map is between two finite sets of the same cardinal (see 1.8(e)) hence it is
a bijection. This proves (f).

\mpb

From (e),(f), we see that we have a bijection $\cf_*(V_{D-1})@>\si>>\cf_*(V_D)$ given by
$E'\m\FF\et_D+E'$ (a refinement of 1.8(e)).

We show:

(g) {\it If $|B|<(D-1)/2$, then there exist $x<y<z$ in $[1,D]-\s$ such that $x,z$ are odd and
$y$ is even.}
\nl
Assume first that any number in $[1,D]-\s$ is even. Then

$1=a_1,b_1,a_2,b_2,\do,a_s,b_s=D$
\nl
are all odd and

$a_2=b_1+2,a_3=b_1+2,\do,a_s=b_{s-1}+2$.
\nl
We have

$b_1=1+c_1, a_2=3+c_1,b_2=3+c_1+c_2,a_3=5+c_1+c_2,b_3=5+c_1+c_2+c_3,\do$,

$a_s=2s-1+c_1+c_2+\do+c_{s-1},b_s=2s-1+c_1+c_2+\do+c_s$.
\nl
The last equality implies $D=2s-1+2|B|-2s$ that is $|B|=(D+1)/2$. This contradicts $|B|<(D-1)/2$.
We see thay the set of odd numbers in $[1,D]-\s$ is nonempty. Let $x$ (resp. $z$) be the
smallest (resp. largest) odd number in $[1,D]-\s$.
We have $b_t<x<a_{t+1}$, $b_u<z<a_{u+1}$ for some $t\in[0,s], u\in[0,s], t\le u$.

Assume now that there is no even number $y\in[1,D]-\s$ such that $x<y<z$. Then

$a_1,b_1,a_2,b_2,\do,a_t,b_t$ are all odd; $a_{t+1},b_{t+1},a_{t+2},b_{t+2},\do,a_u,b_u$ are
all even; $a_{u+1},b_{u+1},\do,a_s,b_s$ are all odd (also $a_1=1$ if $t>0$, $a_1=2$ if $t=0$,
$b_s=D$ if $u<s$, $b_s=D-1$ if $u=s$);

$a_2=b_1+2,a_3=b_2+2,\do,a_t=b_{t-1}+2$; $a_{t+1}=b_t+3$, $a_{t+2}=b_{t+1}+2$,
$a_{t+3}=b_{t+2}+2,\do,a_u=b_{u-1}+2$; $a_{u+1}=b_u+3,a_{u+2}=b_{u+1}+2,\do,a_s=b_{s-1}+2$.
\nl
From this we deduce as above that $b_s=2s-1+c_1+c_2+\do+c_s+e$ where $e=1$ if $x=z$ and $e=2$
if $x<z$. Hence $D=2s-1+2|B|-2s+e$, that is $|B|=(D+1-e)/2$.
If $e=1$ this is a contradiction since $(D+2)/2\n\ZZ$. If $e=2$ we see that
$|B|=(D-1)/2$, contradicting $|B|<(D-1)/2$. This proves (g).

We show:

(h) {\it If $|B|<(D-1)/2$ and $J\in\ci_D$ satisfies $e_J\in E+\FF\et_D$, then $J\in B$.}
\nl
Assume that $J\n B$. By 1.2(f) we then have $e_J\n E$ so that $e_J+\et_D\in E$.
For any $\x\in V_D,i\in[1,D]$ let $(e_i:\x)$ be the coefficient of $e_i$ in $\x$.
Since any element of $E$ is a linear combination of $e_i,i\in\s$, we have
$(e_i:e_J)+(e_i:\et_D)=0$ for $i\n\s$. Let $x<y<z$ be as in (g). Since $x,y,z$ are not in $\s$
we have

$(e_x:e_J)+(e_x:\et_D)=0$, $(e_y:e_J)+(e_y:\et_D)=0$ $(e_z:e_J)+(e_z:\et_D)=0$.
\nl
Since $x,z$ are odd and $y$ is even we have 

$(e_x:\et_D)=1$, $(e_y:\et_D)=0$, $(e_z:\et_D)=1$.
\nl
Hence $(e_x:e_J)=1$, $(e_y:e_J)=0$ $(e_z:e_J)=1$. Thus $e_x$, $e_z$ appear with nonzero
coefficient in $e_J$ so that $x\in J,z\in J$. Since $J$ is an interval and $x<y<z$, it follows
that $y\in J$ contradicting $(e_y:e_J)=0$. This proves (h).

\subhead 2.2\endsubhead
We no longer assume that $D$ is odd.
Let $\ch$ (resp. $\ch'$) be the set of all $i\in[1,s+1]$ such that $a_i-b_{i-1}\ge3$
(resp. $a_i-b_{i-1}\ge4$). We have $\ch'\sub\ch$.

Assume first that $\ch=\emp$ so that $a_i-b_{i-1}=2$ for $i\in[1,s+1]$. From (a) we then have
$$\align&2|B|=2s-1+b_1-a_2+b_2-\do+b_{s-1}-a_s+D\\&=2s-1-2-2-\do-2+D=2s-1-2(s-1)+D=D+1\endalign$$
so that $D$ is odd and $|B|=(D+1)/2$.
We now assume that either $D$ is even or $D$ is odd and $|B|<(D+1)/2$. Then $|\ch|\ge1$.
We write the elements of $\ch$ in a sequence $i_1<i_2<\do<i_t$ with $t\ge1$.
Let $z(B)$ be the set consisting of the
intervals $[a_{i_u}-1,b_{i_{u+1}-1}+1]$ with $u\in[1,t-1]$
and of the intervals $[j,j]$ with $j\in\cup_{i\in\ch'}[b_{i-1}+2,a_i-2]$.
Note that the condition that $z(B)=\emp$ is the same as the condition that
$\ch'=\emp$ and $|\ch|=1$. An equivalent condition is that for some $j\in[1,s+1]$ we have
$a_j-b_{j-1}=3$ and $a_i-b_{i-1}=2$ for all $i\in[1,s+1]-\{j\}$.
From (a) we then have
$$\align&2|B|=2s-2+b_1-a_2+b_2-\do+b_{s-1}-a_s+D=2s-2-2-\do-2+D\\&=2s-2s+D=D\endalign$$
so that $D$ is even and $|B|=D/2$.

Until the end of 2.6 we assume further
that either $D$ is even and $|B|<D/2$ or $D$ is odd and $|B|<(D+1)/2$.
Then $|\ch|+|\ch'|\ge2$ and $z(B)\ne\emp$. We have
$$\align&|z(B)|=\sum_{i\in\ch'}(a_i-b_{i-1}-3)+|\ch|-1\\&=\sum_{i\in\ch}(a_i-b_{i-1}-3)+|\ch|-1
=\sum_{i\in\ch}(a_i-b_{i-1}-2)-1.\endalign$$
If $i\in[1,s+1]-\ch$, then $a_i-b_{i-1}-2=0$, hence the last sum over $\ch$
does not change if $\ch$ is replaced by $[1,s+1]$, so that 
$$\align &|z(B)|=\sum_{i\in[1,s+1]}(a_i-b_{i-1}-2)-1\\&
=-b_0+\sum_{i\in[1,s]}(a_i-b_i)+a_{s+1}-2(s+1)-1\\&
=D-\sum_{i\in[1,s]}(b_i-a_i+2)=D-2|B|.\endalign$$
(We have used 2.1(a).) We set $\D=|z(B)|$, so that $\D=D-2|B|\ge1$.
From the definition we see that there is a well defined sequence
$c_0,c_1,\do,c_\D$ in $\ZZ_{>0}$ such that the intervals in $z(B)$ are:  
$$\align&I_1=[c_0,c_0+c_1-1],I_2=[c_0+c_1,c_0+c_1+c_2-1],\do,\\&
I_\D=[c_0+c_1+\do+c_{\D-1},c_0+c_1+\do+c_{\D}].\tag a\endalign$$  
We write $I_*=(I_1,I_2,\do,I_\D)\in \ci_D\T\ci_D\T\do\T\ci_D$ ($\D$ factors).

\subhead 2.3\endsubhead
We now give several examples of the assignment $B\m I_*$ in 2.2.

When $B=\{[3,5],\{4\},[8,10],\{9\}\}$, then $I_*=(\{1\},[2,6])$ (if $D=10$),
$I_*=(\{1\},[2,6],[7,11])$ (if $D=11$), $I_*=(\{1\},[2,6],[7,11],\{12\})$ (if $D=12$),
$I_*=(\{1\},[2,6],[7,11],\{12\},\{13\})$ (if $D=13$).

If $D=10$ and $B=\{[2,4],\{3\},[8,10],\{9\}\}$ then $I_*=([1,5],\{6\})$.

If $D=10$ and $B=\{[2,4],\{3\},[6,8],\{7\}\}$ then $I_*=([1,9],\{10\})$.

If $D=20$ and $B=\{[4,6],\{5\},[9,11],\{10\},[15,17]\},\{16\}\}$, then

$I_*=(\{1\},\{2\},[3,7],[8,12],\{13\},[14,18],\{19\},\{20\}).$

If $D\ge1,B=\emp$ then $I_*=(\{1\},\{2\},\do,\{D\})$.      

\subhead 2.4\endsubhead
Let $I_*=(I_1,I_2,\do,I_\D)$ be asociated to $B$ as in 2.2. 
From 2.2(a) we see that $\{e_{I_c};c\in[1,\D]\}$ are linearly independent in $V_D$ and
that for $c,d$ in $[1,\D]$:

(a) {\it $(e_{I_c},e_{I_d})$ is $1$ if $c-d=\pm1$ and is $0$ if $c-d\ne\pm1$;
moreover, if $D$ hence $\D$ is odd then $\sum_{c\in[1,\D]}e_{I_c}=e_{[c_0,c_\D]}$
satisfies $(e_{[c_0,c_\D]},e_{I_c})=0$ for all $c\in[1,\D]$.}
\nl
Here $c_0,c_\D$ are as in 2.2(a). Let $\fT_E=\la e_{I_c};c\in[1,\D]\ra\sub V_D$. From (a) we
see that

(b) {\it $\{e_{I_c};c\in[1,\D]\}$ is a basis of $\fT_E$ and $\fT_E^\pe\cap\fT_E$ is $0$
if $D$ is even and is the line $\FF e_{[c_0,c_\D]}$ if $D$ is odd.}
\nl
From the definitions we see that $(e_I,e_{I_c})=0$ for any $I\in B,c\in[1,\D]$.
It follows that $\fT_E\sub E^\pe$. Hence if $x\in E\cap\fT_E$ then $x\in\fT_E^\pe\cap\fT_E$.
Hence $x=0$ if $D$ is even and $x\in\FF e_{[c_0,c_\D]}$ if $D$ is odd.
From the definitions we see that some $c\in[c_0,c_\D]$ is not contained in $\cup_{I\in B}I$;
it follows that $e_{[c_0,c_\D]}$ is not contained in $\EE_B=E$. Since $x\in E$,
it follows that $x=0$. Thus we have $E\cap\fT_E=0$ so that
$E+\fT_E=E\op\fT_E$ has dimension $\dim(E)+\D=|B|+\D=|B|+D-2|B|=D-\dim(E)$.
We have $\dim(E^\pe)=D-\dim(E)$. (When $D$ is even this follows from the nondegeneracy of $(,)$;
when $D$ is odd this follows from 2.1(c) since $|B|\le(D-1)/2$.) We see that
$\dim(E^\pe)=\dim(E\op\fT_E)$.
We show that

(c) $E^\pe=E\op\fT_E$.
\nl
In view of the dimension equality above it is enough to show that $E\op\fT_E\sub E^\pe$.
The inclusion $E\sub E^\pe$ holds since $E$ is isotropic (see 1.9). To prove the
inclusion $\fT_E\sub E^\pe$ it is enough to prove that $(e_{I_c},e_I)=0$ for any
$c\in[1,\D]$ and any $I\in B$. But from 2.2(a) and the definitions we see that
$I_c\prec I$ or $I\prec I_c$ or $I,I_c$ are non-touching. In each case we have 
$(e_{I_c},e_I)=0$, proving (c).

We see that $\fT_E$ is a canonical complement of $E$ in $E^\pe$ and that $\fT_E$
has a canonical basis $\{e_{I_c},c\in[1,\D]\}$.

\subhead 2.5\endsubhead
We now assume that $D\ge2$ and $i\in[1,D]$, $B'\in S_{D-2}$
are such that $B=t_i(B')$. Let $E'\in\cf(V_{D-2})$ be the subspace defined by $B'$ so that
$E=T_i(E')\op\FF e_i$. We have $|B'|=|B|-1$ hence
$|B'|<(D-2)/2$ if $D$ is even and $|B'|<(D-1)/2$ if $D$ is odd. Thus
$\fT_{E'}\sub V_{D-2}$ is defined.
From the definitions we see that

(a) {\it $T_i:V_{D-2}@>>>V_D$ carries $\fT_{E'}$ isomorphically
onto $\fT_E$ compatibly with the canonical bases of $\fT_{E'}$ and $\fT_E$.}

\subhead 2.6\endsubhead
For $c\le d$ in $[1,\D]$ we set
$$I_{[c,d]}=I_c\cup I_{c+1}\cup I_{c+2}\cup\do\cup I_d\in\ci_D.$$
If $\D$ is even and $k\in[1,\D/2]$ we set
$$\fT^k_E=\la e_{I_{[1,\D]}},e_{I_{[2,\D-1]}},\do,e_{I_{[k,\D+1-k]}}\ra.$$
If $\D$ is odd and $k\in[1,(\D-1)/2]^1$ let
$$\fT^k_E=\la e_{I_{[1,\D-1]}},e_{I_{[2,\D-2]}},\do,e_{I_{[k,\D-k]}}\ra,$$
$$\fT^{-k}_E=\la e_{I_{[2,\D]}},e_{I_{[3,\D-1]}},\do,e_{I_{[k+1,\D+1-k]}}\ra.$$
Without restriction on $\D$ we set $\fT^0_E=0$.
Thus $\fT^k_E$ is defined for any

(a) $k\in[0,\D/2]$ if $\D$ is even and for
any $k\in[1,(\D-1)/2]^1\cup(-[1,(\D-1)/2]^1)\cup\{0\}$ if $\D$ is odd.
\nl
For $k$ as in (a) we define $E(k)=E+\fT^k_E=E\op\fT^k_E\sub V_D$. We define $B(k)$ to be

$B$ if $k=0$,

$B\sqc\{I_{[1,\D]},I_{[2,\D-1]},\do,I_{[k,\D+1-k]}\}$ if $\D$ is even, $k\in[1,\D/2]$,

$B\sqc\{I_{[1,\D-1]},I_{[2,\D-2]},\do,I_{[k,\D-k]}\}$ if $\D$ is odd, $k\in[1,(\D-1)/2]^1$,

$B\sqc\{I_{[2,\D]},I_{[3,\D-1]},\do,I_{[-k+1,\D+1+k]}\}$ if $\D$ is odd, $-k\in[1,(\D-1)/2]^1$.
\nl
We show:

(b) {\it $B(k)\in\SS_D$ and $E(k)=\EE_{B(k)}\in\pmb\cf(V_D)$}.
\nl
We argue by induction on $D$. If $E=0$ then $E(k)$ is a primitive subspace of $V_D$
so that (b) holds. Thus we can assume that $D\ge2$ and $E\ne0$.
(If $D=1$ then $E=0$ since $|B|<(D+1)/2$.)
We can find $i,B'\in S_{D-2},E'\in\cf(V_{D-2}]$ as in 2.5. By the
induction hypothesis we have $B'(k)\in\SS_{D-2}$ and
$E'(k)\in\pmb\cf(V_{D-2})$ is the subspace of $V_{D-2}$ defined in terms of $B'(k)$
in the same way as $E(k)$ is defined in terms of $B(k)$.
We have $E=T_i(E')\op\FF e_i$. Using 2.5 we see that $E(k)=T_i(E'(k))\op\FF e_i$
and $B(k)=t_i(B'(k))$. Hence (b) holds.

The same inductive proof shows that

(c) {\it if $I\in B(k)-B$ then $|I|$ is even.}
\nl
(We use that this holds when $B=\emp$.)

\subhead 2.7\endsubhead
We no longer fix $B$. Let $\un\SS_D$ be the set of all pairs $(B,k)$ where $B\in S_D$ and
one of (i)-(iii) below holds.

(i) $|B|=D/2$ and $k=0$ (if $D$ is even);

(ii) $|B|=(D+1)/2$ and $k=0$ (if $D$ is odd);

(iii) we have $|B|<D/2$ (if $D$ is even), $|B|<(D+1)/2$ (if $D$ is odd)
so that $\D:=D-2|B|\ge1$, and we have $k\in[0,\D/2]$ if $D$ is even
and $k\in[1,(\D-1)/2]^1\cup(-[1,(\D-1)/2]^1)\cup\{0\}$ if $D$ is odd.
\nl
We define $\l:\un\SS_D@>>>\SS_D$ by $\l(B,k)=B(k)$. (We have $B(0)=B$.)

\proclaim{Proposition 2.8} $\l$ is a bijection.
\endproclaim
We show that $\l$ is injective. Indeed, assume that $(B,k)\in\un\SS_D$,
$(B',k')\in\un\SS_D$ are such that $B(k)=B'(k')$. Intersecting the two sides with $\ci_D^1$
and using 2.6(c) we obtain $B=B'$. Since $|B(k)|-|B|=\pm k$, $|B(k')|-|B|=\pm k'$, we see
that $\pm k=\pm k'$. If $D$ is even we have $k\ge0,k'\ge0$ hence $k=k'$. Assume now that
$D$ is odd and $k\ne k'$ so that $k'=-k$. But then one of $B(k)-B$, $B(k')-B$
contains $I_{[1,\D-1]}$ while the other one does not; this contradicts
$B(k)-B=B(k')-B$. We see that we must have $k=k'$ proving the injectivity of $\l$.

We show that $\l$ is surjective by induction on $D$.
Let $B\in\SS_D$. If $B$ is primitive then $B=\{0\}(k)$ with $k\in[0,D/2]$ (if $D$ is even)
or with $k\in[1,(D-1)/2]^1\cup(-[0,(D+1)/2]^1)\cup\{0\}$ (if $D$ is odd); thus
$B$ is in the image of $\l$. If $B\in S_D$ then $B=B(0)$ so that $B$ is in the image of $\l$.
Now assume that $B$ is not primitive and $B\n S_D$.
Then $D\ge2$ and there exists $i\in[1,D]$ and $B'\in\SS_{D-2}$ such that $B=t_i(B')$.
We must have $B'\n S_{D-2}$.
By the induction hypothesis we have
$B'=\tB'(k)$ where $\tB'\in S_{D-2}$ and
$k\in[1,(D-2)/2-|\tB'|]$ (if $D$ is even) and
$k\in[1,(D-2-2|\tB'|-1)/2]^1\cup(-[1,(D-2-2|\tB'|-1)/2]^1)$   (if $D$ is odd).
Let $\tB=t_i(\tB')$. We have $\tB\in S_D$ and $|\tB|=|\tB'|+1$ so that
$k\in[1,D/2-|\tB|]$ (if $D$ is even) and
$k\in[1,(D-2|\tB|-1)/2]^1\cup(-[1,(D-2|\tB|-1)/2]^1)$   (if $D$ is odd).
As in the proof of 2.6(b) we have $\tB(k)=t_i(\tB'(k))$ hence $\tB(k)=B$. Thus $\l$
is surjective. The proposition is proved.

\subhead 2.9\endsubhead
In this subsection we assume that $D$ is odd. Let
$$S'_D=S_{D-1}\sqc\{B\in S_D-S_{D-1};|B|<(D-1)/2\}.$$
This is a subset of $S_D$ (we view $S_{D-1}$ as a subset of $S_D$, see 1.10). Let
$$\align&\cf'(V_D)=\{\EE_B;B\in S'_D\}\\&
=\cf(V_{D-1})\sqc\{E\in\cf(V_D)-\cf(V_{D-1});\dim(E)<(D-1)/2\}.\endalign$$
This is a subset of $\cf(V_D)$ (we view $\cf(V_{D-1})$ as a subset of $\cf(V_D)$, see 1.10).

Let $\SS'_D$ be the set consisting of all elements $B(k)\in\SS_D$ (see 2.6) where either

(i) $B\in S_{D-1}$ and $k\in\{0\}\cup[1,(D-2|B|-1)/2]^1$ or

(ii) $B\in S_D-S_{D-1}$, $|B|<(D-1)/2$ and $k\in[1,(D-2|B|-1)/2]^1$.
\nl
(We use that $S_{D-1}\sub S_D$, see 1.10.)
We define a map $\io:\SS'_D@>>>\SS_{D-1}$ by $B(k)\m B'(k')$ where $B'(k')\in\SS_{D-1}$
(as in 2.6 with $D$ replaced by $D-1$) is given by:

$B'=B,k'=0$ if $B,k$ are as in (i) and $k=0$;

$B'=B,k'=(k+1)/2$ if $B,k$ are as in (i) and $k>0$;

$B'\in S_{D-1}$ is obtained by removing from $B$ the unique interval containing $D$,
$k'=(k+3)/2$ if $B,k$ are as in (ii).
\nl
From the definitions we see using 2.8 that

(a) {\it $\io$ is a bijection.}

\subhead 2.10\endsubhead
Let $\un{\pmb\cf}(V_D)$ be the set of all pairs $(E,k)$ where $E\in\cf(V_D)$ and
one of (i)-(iii) below holds.

(i) $\dim(E)=D/2$ and $k=0$ (if $D$ is even);

(ii) $\dim(E)=(D+1)/2$ and $k=0$ (if $D$ is odd);

(iii) we have $\dim(E)<D/2$ (if $D$ is even), $\dim(E)<(D+1)/2$ (if $D$ is odd)
so that $\D:=D-2\dim(E)\ge1$, and we have $k\in[0,\D/2]$ if $D$ is even
and $k\in[1,(\D-1)/2]^1\cup(-[1,(\D-1)/2]^1)\cup\{0\}$ if $D$ is odd.
\nl
We define $\mu:\un{\pmb\cf}(V_D)@>>>\pmb\cf(V_D)$
by $\mu(E,k)=E(k)$. (We have $E(0)=E$.)
The following result is a reformulation of 2.8.

(a) {\it $\mu$ is a bijection.}

\subhead 2.11\endsubhead
In the remainder of this section we assume that $D$ is odd. We set
$$V'_D=V_D/\FF\et_D$$
where $\et_D\in V_D$ is as in 1.8; let $\p:V_D@>>>V'_D$ be the obvious map.
Let $\cf(V'_D)$ be the set of subspaces of $V'_D$ of the form $\p(E)$ for various $E\in\cf(V_D)$.
We show:

(a) {\it the map $E\m\p(E)$ defines a bijection $\cf'(V_D)@>\si>>\cf(V'_D)$ (notation of 2.9).}
\nl
We first show that this map is surjective. It is enough to show that for any $E\in\cf(V_D)$
there exists $E'\in\cf'(V_D)$ such that $\p(E)=\p(E')$.
If $E\in\cf'(V_D)$ then $E'=E$ satisfies our requirement. Thus we can assume that
$E\n\cf'(V_D)$ so that $E\in\cf(V_D)-\cf(V_{D-1})$ and $\dim(E)\ge(D-1)/2$.
It follows that we have either $\dim(E)=(D+1)/2$ or $\dim(E)=(D-1)/2$.
If $\dim(E)=(D+1)/2$ then by 2.1(f) we can find $E'\in\cf(V_{D-1})$ such that $\p(E)=\p(E')$;
since $E'\in\cf'(V_D)$ we see that $E'$ satisfies our requirement.
If $\dim(E)=(D-1)/2$ then by 2.1(e) we have $\et_D\n E$ and $\FF\et_D\op E\in\cf(V_D)$; by the
previous sentence applied to $\FF\et_D\op E$ (which has dimension $(D+1)/2$) instead of $E$
we see that we can find $E'\in\cf(V_{D-1})$ such that $\p(\FF\et_D+E)=\p(E')$. Since
$\p(\FF\et_D+E)=\p(E)$ we see that $E'$ satisfies our requirement.
This proves the surjectivity of the map in (a).

We now prove injectivity. Let $E,E'$ in $\cf'(V_D)$ be such that $\p(E)=\p(E')$; we must
show that $E=E'$.

Let $B\in S'_D,B'\in S'_D$ be such that $E=\EE_B,E'=\EE_{B'}$. Since $\p(E)=\p(E')$ we have
$E+\FF\et_D=E'+\FF\et_D$. Assume first that $\dim(E')<(D-1)/2$. 
If $J\in B$, then $e_J\in E$ hence $e_J\in E+\FF\et_D=E'+\FF\et_D$.
Using 2.1(h), we deduce that $J\in B'$. Thus we have $B\sub B'$ so that $E\sub E'$ and
$\dim(E)<(D-1)/2$. Similarly, if $\dim(E)<(D-1)/2$, then $E'\sub E$ and $\dim(E')<(D-1)/2$. We
see that if one of the conditions $\dim(E)<(D-1)/2$, $\dim(E')<(D-1)/2$ holds, then both
conditions hold and $E=E'$. Thus, we can assume that $E\in\cf(V_{D-1}),E'\in\cf(V_{D-1})$.
In this case we use that the restriction of $\p$ to $V_{D-1}$ is an isomorphism
$V_{D-1}@>\si>>V'_D$ to see that $E=E'$. This proves the injectivity of the map in (a) and
completes the proof of (a).

\mpb

We have $V'_D=V'_D{}^0\op V'_D{}^1$ where $V'_D{}^0$ (resp. $V'_D{}^1$) is the image of
$V_D^0$ (resp. $V_D^1$) under the map $\p:V_D@>>>V'_D$.

Let $\cc(V'_D{}^1)$ be the set of subspaces of $V'_D{}^1$ which are images under $\p$ of
subspaces of $V_D^1$ in $\cc(V_D^1)$.

Let $\ovm{\cc}(V'_D{}^1)$ be the set of pairs of subspaces of $V'_D{}^1$ which are images under
$\p$ of pairs of subspaces of $V_D^1$ in $\ovm{\cc}(V_D^1)$.
We have $\ovm{\cc}(V'_D{}^1)\sub\cc(V'_D{}^1)\T\cc(V'_D{}^1)$.

By 1.9(b), $\ovm{\cc}(V'_D{}^1)$ consists of the pairs $\g_E=(\p(E^1),\p((E^0)^!))$
of subspaces of $V'_D{}^1$ for various $E\in\cf(V_D)$. We show: 

(b) {\it for $E,E'$ in $\cf(V_D)$ we have $\g_E=\g_{E'}$ if and only if $\p(E)=\p(E')$.}
\nl
Assume first that $\g_E=\g_{E'}$ that is $\p(E^1)=\p(E'{}^1)$ and $\p((E^0)^!)=\p((E'{}^0)^!)$.
We have $E^1+\FF\et_D=E'{}^1+\FF\et_D$, $(E^0)^!=(E'{}^0)^!$ so that
$E^0=E'{}^0$ and $E+\FF\et_D=E'+\FF\et_D$ that is $\p(E)=\p(E')$.
Conversely assume that $\p(E)=\p(E')$. Then $E+\FF\et_D=E'+\FF\et_D$ hence
$E^1+\FF\et_D=E'{}^1+\FF\et_D$, $E^0=E'{}^0$, so that $(E^0)^!=(E'{}^0)^!$ and $\g_E=\g_{E'}$.
This proves (b).

\mpb

From (b) we see that $E\m\g_E$ induces a bijection $\cf(V'_D)@>\si>>\ovm{\cc}(V'_D{}^1)$.
Composing this bijection with the bijection (a), we obtain a bijection

(c) $\cf'(V_D)@>\si>>\ovm{\cc}(V'_D{}^1)$. 

\subhead 2.12\endsubhead
Let $\pmb\cf'(V_D)=\{\EE_B;B\in\SS'_D\}\sub\pmb\cf(V_D)$.
Let $\un{\pmb\cf'}(V_D)$ be the set of all pairs $(E,k)$ where 

(i) $E\in\cf(V_{D-1})$ and $k\in\{0\}\cup[1,(D-2\dim(E)-1)/2]^1$ or

(ii) $E\in\cf(V_D)-\cf(V_{D-1})$, $\dim(E)<(D-1)/2$ and $k\in[1,(D-2\dim(E)-1)/2]^1$.
\nl
We define $\mu':\un{\pmb\cf'}(V_D)@>>>\pmb\cf'(V_D)$ by $\mu'(E,k)=E(k)$. The following result
is a consequence of 2.10(a).

(a) {\it $\mu'$ is a bijection.}

\subhead 2.13\endsubhead
Let $\cf_*(V'_D)=\{E'\in\cf(V'_D);\dim(E')=(D-1)/2\}$. We show:

(a) {\it The map $E\m\p(E)$ defines a bijection $\cf_*(V_{D-1})@>>>\cf_*(V'_D)$.}
\nl
Let $E\in\cf_*(V_{D-1})$. Since the restriction of $\p$ to $V_{D-1}$ is bijective we see that
$\dim(\p(E))=(D-1)/2$. Using 2.11(a) we see that the map in (a) is well defined and injective.
Let $E'\in\cf_*(V'_D)$. By 2.11(a) we have $E'=\p(E)$ where $E\in\cf'(V_D)$. We have
$\dim(E)\ge(D-1)/2$ hence $E\in\cf(V_{D-1})$ and $\dim(E)=(D-1)/2$ so that $E\in\cf_*(V_{D-1})$.
We see that the map in (a) is surjective. This proves (a).

\mpb

The following result can be deduced from (a) and 2.1(e),(f).

(b) {\it The map $E\m\p(E)$ defines a bijection $\cf_*(V_D)@>>>\cf_*(V'_D)$.}

\head 3. From  $\SS_D$ to $V_D$\endhead
\subhead 3.1\endsubhead                    
We fix a symbol $\aa$. We set $\ov{[1,D]}=[1,D]\sqc\{\aa\}$.
A subset $I$ of $\ov{[1,D]}$ is said to be an interval if
either $I\in\ci_D$ or if $\ov{[1,D]}-I\in\ci_D$.
Such an $I$ is necessarily nonempty and not equal to $\ov{[1,D]}$.
Let $\ti\ci_D$ be the set of intervals of $\ov{[1,D]}$.
We have $\ci_D\sub\ti\ci_D$.
For $I\in\ci_D$ we define $I^*\in\ti\ci_D$ by $I^*=I$ if $|I|$ is odd and $I^*=\ov{[1,D]}-I$
if $|I|$ is even.
For $B\in R_D$ we set $B^*=\{I^*;I\in B\}$.
Let $\SS_D^*=\{B^*;B\in\SS_D\}$.
If $D$ is even and $B\in\SS_D^*$ then any $I\in B$ satisfies $|I|=1\mod2$. 
If $D$ is odd and $B\in\SS_D^*$ then any $I\in B$ satisfies $|I|=1\mod2$ (if $\aa\n I$) or
$|I|=0\mod2$ (if $\aa\in I$).

\subhead 3.2\endsubhead
Let $B\in\SS_D^*$. For any $j\in\ov{[1,D]}$ we set $n_j(B)=\sha(I\in B;j\in I)$.
We now define $\e^j(B)\in\FF$ for any $j\in\ov{[1,D]}$ as follows.
If $D$ is even we set $\e^j(B)=n_j(B)(n_j(B)+1)/2\in\FF$.
Assume now that $D$ is odd. We denote by $[B]$ the union of all $I\in B$ such that $\aa\in I$;
then either $[B]=\emp$ or $\aa\in[B]$ and $[B]\in B$. If $j\n[B]$ or if $j=\aa$ we set
$\e^j(B)=n_j(B)(n_j(B)+1)/2\in\FF$. Assume now that $j\in[B],j\ne\aa$.
We have $[B]-\{\aa\}=[1,k]\sqc[k',D]$ 
with $0\le k<D$, $1<k'\le D+1$, $k'\ge k+2$; we have either $j\in[k',D]$ (we then set $u=k'$) or   
$j\in[1,k]$ (we then set $u=k$).
If $u=j+1\mod2$ we set $\e^j(B)=(n_j(B))(n_j(B)+1)/2\in\FF$.
If $u=j\mod2$ we set $\e^j(B)=(n_j(B)+1)(n_j(B)+2)/2\in\FF$.

We set $e_\aa=e_{[1,D]}\in V_D$; we define $\e_D:\SS_D^*@>>>V_D$ by
$\e_D(B)=\sum_{j\in\ov{[1,D]}}\e^j(B)e_j$.

In the case where $D$ is odd, we define $\e'_D:\SS'{}^*_D@>>>V'_D$ by $B\m\p(\e_D(B))$. We state
the following result (see also the tables in \S4).

\proclaim{Theorem 3.3}(a) For any $D$, $\e_D:\SS^*_D@>>>V_D$ is a bijection.

(b) For $D$ odd, $\e'_D:\SS'{}^*_D@>>>V'_D$ is a bijection.
\endproclaim
When $D$ is even, (a) can be deduced from \cite{L20, 1.17(b)}, see 3.4; the proof of (a)
for odd $D$ is similar and will be omitted.

We prove (b). Consider the diagram
 $$\CD
 \SS'_D@>\io>>\SS_{D-1}\\
@VVV         @VVV   \\
V_D @<<< V_{D-1}
\endCD$$
with $\io$ being the bijection in 2.9, the left vertical map being $B\m\e_D(B^*)$,
the right vertical map being $B\m\e_{D-1}(B^*)$,
and the lower horizontal map being the obvious inclusion. From the definitions this
diagram is commutative.
Combining this with (a) (with $D$ replaced by $D-1$) and with the bijection
$\SS_D@>>>\SS^*_D$, $B\m B^*$, we see that the restriction of $\e_D$ to $\SS'{}^*_D$
is a bijection $\SS'{}^*_D@>>>V_{D-1}$. It remains to use that $\p$ restricts to a
bijection $V_{D-1}@>>>V'_D$.

\subhead 3.4\endsubhead
In this subsection we assume that $D$ is even. Let $B\in\SS_D$. For $i\in[1,D]$, $\d\in\{0,1\}$,
let $B^\d_i=\{I\in B^\d;i\in I\}$, $\k(B)=|B^0|$,
$\ti\e^i(B)=(|B^1_i|-|B^0_i|-\un{|B^0|})(|B^1_i|-|B^0_i|-\un{|B^0|}+1)/2\in\FF$,
where for $z\in\ZZ$ we define $\un z\in\{0,1\}$ by $z=\un z\mod2$.
We define $\ti\e(B)=\sum_{i\in[1,D]}\ti\e^i(B)e_i\in V_D$.
By \cite{L20, 1.17(b)}, $B\m\ti\e(B)$ defines a bijection $\SS_D@>\si>>V_D$.
Hence to prove 3.3(a) in our case it is enough to show that for any $B\in\SS_D$
we have $\ti\e(B)=\e_D(B^*)$ that is $\ti\e^i(B)=\e^i(B^*)+\e^\aa(B^*)$
for any $i\in[1,D]$, or equivalently that
$$(|B^1_i|-|B^0_i|-\un{|B^0|})(|B^1_i|-|B^0_i|-\un{|B^0|}+1)/2
=n_i(n_i+1)/2+n_\aa(n_\aa+1)/2 \mod2$$
where we write $n_i,n_\aa$ instead of $n_i(B^*),n_\aa(B^*)$.
From the definition we have $n_\aa=|B^0|$,
$$n_i=|B^1_i|+\sha(I\in B^0,i\n I)=|B^1_i|+|B^0|-|B^0_i|.$$
Hence
$$n_i(n_i+1)+n_\aa(n_\aa+1)=(|B^1_i|+|B^0|-|B^0_i|)(|B^1_i|+|B^0|-|B^0_i|+1)+|B^0|(|B^0|+1).$$
Setting $b=|B^0|,p=|B^1_i|-|B^0_i|-\un b$, we see that it is enough to show
$$(p+b+\un b)(p+b+\un b+1)+b(b+1)=p(p+1)\mod4.$$
This follows from $2p(b+\un b)=0\mod4$, $2b^2+2b=0\mod4$, $2b\un b+(\un b)^2+\un b=0\mod4$.
This completes the deduction of 3.3(a) (for $D$ even) from \cite{L20}.

\subhead 3.5\endsubhead
For $I\in\ti\ci_D$ we set $e_I=\sum_{j\in I}e_j\in V_D$ (with $e_a$ as in 3.2.)
Whem $I\in\ci_D$ this agrees with the earlier definition.
For $B\in\SS^*_D$ let $\EE_B$ be the subspace of $V_D$ spanned by $\{e_I;I\in B\}$.
For $B\in\SS_D$ we have

$\EE_B=\EE_{B^*}$.
\nl
(It is enough to show that for any $I\in\ci_D$ we have $e_{I^*}=e_I$; this is clear from the
definition.)

We define $\bar\e_D:\pmb\cf(V_D)@>>>V_D$ by $\EE_B=\EE_{B^*}\m\e_D(B)$ for any $B\in\SS^*_D$.
The following result is a reformulation of 3.3(a) (we use 1.2(g)).

(a) {\it For any $D$, the map $\bar\e_D:\pmb\cf(V_D)@>>>V_D$ is a bijection.}
\nl
We have the following result.

(b) {\it If $\EE\in\pmb\cf(V_D)$, then $\bar\e_D(\EE)\in\EE$; this property characterizes the
bijection $\bar\e_D$.}
\nl
When $D$ is even this can be deduced from \cite{L20, 1.22(b)} using the arguments in 3.4.
A similar proof applies when $D$ is odd.

When $D$ is odd, we
define $\bar\e'_D:\pmb\cf'(V_D)@>>>V'_D$ by $\EE_B\m\e'_D(B)$ for any $B\in\SS'{}^*_D$.
The following result is a reformulation of 3.3(b).

(c) {\it For odd $D$, the map $\bar\e'_D:\pmb\cf'(V_D)@>>>V'_D$ is a bijection.}
\nl
From (b) we deduce (for $D$ odd):

(d) {\it If $\EE\in\pmb\cf'(V_D)$, then $\bar\e'_D(\EE)\in\p(\EE)$.}
\nl
One can show that

(e) {\it the property in (d) characterizes the bijection $\bar\e'_D$.}
\nl
From (e) we can deduce that (for odd $D$):

(f) {\it if $\EE,\EE'$ in $\pmb\cf'(V_D)$ satisfy $\p(\EE)=\p(\EE')$, then $\EE=\EE'$.}
\nl
Indeed, let $\EE,\EE'$ in $\pmb\cf'(V_D)$ be such that  $\p(\EE)=\p(\EE')$.
We define a new bijection $\e'':\pmb\cf'(V_D)@>>>V'_D$ 
by $\e''(\EE_1)=\bar\e'_D(\EE_1)$ if $\EE_1\in\pmb\cf'(V_D)$, $\EE_1\n\{\EE,\EE'\}$,
$\e''(\EE)=\bar\e'_D(\EE')$, $\e''(\EE')=\bar\e'_D(\EE)$.
For $\EE_1\in\pmb\cf'(V_D)$ we have $\e''(\EE_1)\in\p(\EE_1)$. (When $\EE_1\n\{\EE,\EE'\}$,
this follows from (d); when $\EE_1\in\{\EE,\EE'\}$, this follows (d) and from
$\p(\EE)=\p(\EE')$.) Using (e) we deduce that $\e''=\bar\e'_D$ so that
$\bar\e'_D(\EE)=\bar\e'_D(\EE')$. Using the injectivity of $\bar\e'_D$ we deduce that
$\EE=\EE'$, as desired.

\subhead 3.6\endsubhead
Note that

(a) {\it there exists a partial order $\le$ on $V_D$ such that for any $B\in\SS^*_D$ and any
$x\in\EE_B$ we have $x\le\e_D(B)$.}
\nl
When $D$ is even this follows from \cite{L20,1.22(b)}. The proof for $D$ odd is similar.

Assuming that $D$ is odd, one can show that the following analogue of (a) holds:

(b) {\it there exists a partial order $\le$ on $V'_D$ such that for any $B\in\SS'{}^*_D$ and any
$x\in\p(\EE_B)$ we have $x\le\e'_D(B)$.}

\subhead 3.7\endsubhead
We now assume that $G$ is as in 0.6(i). We prove 0.6(a) in
this case. By \cite{L84} we can identify $\cg_c$ with $V_{D'}^\d$ for some even $D'\ge0$ and some
$\d\in\{0,1\}$. We can assume that $D'=D$. Then $M(\cg_c)$ becomes
$$V_D^\d\op\Hom(V_D^\d,\FF)=V_D^{1-\d}\op V_D^\d=V_D.$$
By \cite{L82, L86, L87}, we have $\HH_c=\cc(V_D^d)$. We define $\ovm{\HH}_c$ to be the subset
$\ovm{\cc}(V_D^\d)$ of $\cc(V_D^\d)\T\cc(V_D^\d)=\HH_c\T\HH_c$.
Let $(\cl,\tcl)\in\ovm{\cc}(V_D^\d)=\ovm{\HH}_c$ (we have $\cl\sub\tcl$).
Let $E=\cl\op\tcl^!\sub V_D$. We have $E\sub E^\pe=\tcl\op\cl^!$. Recall that $E\in\cf(V_D)$.
We have
$$M(\tcl/\cl)=\tcl/\cl\op\Hom(\tcl/\cl,\FF)=\tcl/\cl\op\cl^!/\tcl^!=
(\tcl\op\cl^!)/(\cl\op\tcl^!)=E^\pe/E.$$
Let $Prim(\cl,\tcl)=Prim(E^\pe/E)$ be the subset of $M(\tcl/\cl)=E^\pe/E$ consisting of the
subspaces $E(k)/E$ of $E^\pe/E$ for various $k\in[0,D/2-\dim(E)]$ (notation of 2.6). Since
$(\cl,\tcl)\m E$ identifies $\ovm{\cc}(V_D^\d)$ with $\cf(V_D)$, we see that the map
appearing in 0.6(a) can be identified with a map $\op_{E\in\cf(V_D)}Prim(E^\pe/E)@>>>V_D$.
We define such a map by $(E,E(k)/E)\m\bar\e_D(E(k))$.
Using 2.10(a) and 3.3 we see that this map is a bijection. From the definitions, if
$E$ corresponds to $(\cl,\tcl)$ as above, then $\ss_{\cl,\tcl;V_D^\d}(E(k)/E)\in\CC[V_D]$ is the
characteristic function of $E(k)$ hence its support contains $\bar\e_D(E(k))$ (see 3.5(b)).
This proves the existence of the bijection $\Th$ in 0.6(a) in our case. The uniqueness of $\Th$
follows from 3.5(b).

\subhead 3.8\endsubhead
We now assume that $G$ is as in 0.6(ii). We prove 0.6(a) in this case. By \cite{L84} we can
identify $\cg_c$ with $V'_{D'}{}^1$ (see 2.11) for some odd $D'\ge1$. We can assume that
$D'=D$. Then $M(\cg_c)$ becomes
$$V'_D{}^1\op\Hom(V'_D{}^1,\FF)=V'_D{}^0\op V'_D{}^1=V'_D.$$

By \cite{L82, L86, L87}, we have $\HH_c=\cc(V'_D{}^1)$ (see 2.11).
We define $\ovm{\HH}_c$ to be the subset $\ovm{\cc}(V'_D{}^1)$ of
$\cc(V'_D{}^1)\T\cc(V'_D{}^1)=\HH_c\T\HH_c$ (see 2.11).

Consider an element of $\ovm{\cc}(V'_D{}^1)=\ovm{\HH}_c$; it can be written in the
form $\g_E=(\p(E^1),\p((E^0)^!)$ where $E\in\cf'(V_D)$ is uniquely determined (see 2.11). We have
$E\sub E^\pe=(E^0)^!\op(E^1)^!$ and
$$\align&M(\p((E^0)^!)/\p(E^1))=\p((E^0)^!)/\p(E^1)\op\Hom(\p((E^0)^!)/\p(E^1),\FF)\\&
=\p((E^0)^!)/\p(E^1)\op\p((E^1)^!)/\p(E^0)\\&=(\p((E^0)^!)\op\p((E^1)^!))/(\p(E^1)\op\p(E^0))=
\p(E^\pe)/\p(E).\endalign$$

Let $Prim(\g_E)$ be the subset of $M(\p((E^0)^!)/\p(E^1))=\p(E^\pe)/\p(E)$
consisting of the subspaces $\p(E(k))/\p(E)$ (notation of 2.6) of $\p(E^\pe)/p(E)$
for various $k$ in 

$[1,(D-2\dim(E)-1)/2]^1$ if $E\in \cf'(V_D)-\cf(V_{D-1})$, $\dim(E)<(D-1)/2$, or in

$\{0\}\cup[1,(D-2\dim(E)-1)/2]^1$ if $E\in\cf(V_{D-1})$.
\nl
We now see that the map appearing in 0.6(a) can be identified with a map
$\op_{E\in\cf'(V_D)}Prim(\g_E)@>>>V'_D$. We define such a map by $(E,E(k)/E)\m\bar\e'_D(E(k))$.
Using 2.12(a) and 3.5(c) we see that this map is a bijection. From the definitions,
for $E,k$ as above we see that $\ss_{\g_E;V'_D{}^1}(\p(E(k))/\p(E))\in\CC[V'_D]$ is the
characteristic function of $\p(E(k))$ hence its support contains
$\bar\e'_D(\p(E(k)))$ (see 3.5(d)). This proves the existence of the bijection $\Th$ in 0.6(a)
(in our case). The uniqueness of $\Th$ follows from 3.5(e).

\subhead 3.9\endsubhead
Following \cite{L20, 1.11} we define $u:V_D@>>>\ZZ$ by
$$u(x)=\sum_{s\in[1,r];a_s\ne b_s\mod2}(-1)^{a_s}$$
where
$$x=e_{[a_1,b_1]}+e_{[a_2,b_2]}+\do+e_{[a_r,b_r]}\in V_D$$
with $1\le a_1\le b_1\ll<a_2\le b_2\ll\do\ll a_r\le b_r\le D$.   
When $D$ is even let $\tu:V_D@>>>\NN$ be the function defined
by $\tu(x)=2u(x)$ if $u(x)\ge0$, $\tu(x)=-2u(x)-1$ if $u(x)<0$.
When $D$ is odd let $\tu:V_D@>>>\NN$ be the function defined
by $\tu(x)=2u(x)-1$ if $u(x)>0$, $\tu(x)=-2u(x)+1$ if $u(x)<0$, $\tu(x)=0$ if $u(x)=0$.

(a) {\it Let $k\in\NN$. Under the bijection $\e_D:\SS^*_D@>>>V_D$ (see 3.3(a)),
the subset $\tu\i(k)$ of $V_D$ corresponds to the subset $\SS^*_D(k)$ of $\SS^*_D$
consisting of all $B\in\SS^*_D$ such that the number of intervals $I\in B$
containing $\aa$ (or equivalently such that $|I^*|$ is even) is equal to $k$.}
\nl
When $D$ is even, this follows from \cite{L20, 1.14(c)}; the case where $D$ is odd is similar.

\subhead 3.10\endsubhead
In this subsection we assume that $D$ is odd. We define $\tu':V'_D@>>>\NN$ as the composition
$V'_D@>>>V_{D-1}@>>>V_D@>>>\NN$ where the
first map is the inverse of the bijection $V_{D-1}@>>>V'_D$ induced by $\p$, the second map is
the obvious inclusion and the third map is $\tu$ as in 3.9.
The following result can be deduced from 3.9(a).

(a) {\it Let $k\in\NN$. Under the bijection $\e'_D:\SS'{}^*_D@>>>V'_D$ (see 3.3(b)), the subset
$\tu'{}\i(k)$ of $V'_D$ corresponds to the subset $\SS'{}^*_D(k)=\SS'{}^*_D\cap\SS^*_D(k)$ of
$\SS'{}^*_D$ (see 3.9(a)).}

\subhead 3.11\endsubhead
In this subsection we assume that $D$ is even.
We define a bijection $\k:\ov{[1,D]}@>>>\ov{[1,D]}$ by $\k(j)=j+1$ if $j\in[1,D-1]$,
$\k(D)=\aa$, $\k(\aa)=1$. This induces a bijection $\k:\ti\ci_D@>>>\ti\ci_D$,
$I\m\k(I)=\{\k(j);j\in I\}$ and a bijection $\k:R_D@>>>R_D$, $B\m\k(B)=\{\k(I);I\in B\}$.
From \cite{L20a} it is known that for $B\in R_D$ we have

(a) {\it $B\in\SS^*_D$ if and only if $\k(B)\in\SS^*_D$.}
\nl
We show:

(b) {\it Let $B\in R_D$ be such that $|I|$ is odd for all $I\in B$. We have $B\in\SS^*_D$ if and
only if $\k^s(B)\in S_D$ for some $s\ge0$.}
\nl
We can assume that $D\ge2$. 
If $B$ satisfies $\k^s(B)\in S_D$ for some $s\ge0$ then by (a) we have $B\in\SS^*_D$ (since
$S_D\sub\SS^*_D$). Conversely, assume that $B\in\SS^*_D$. From the definition of $\SS^*_D$
we see by induction on $D$ that $\ov{[1,D]}-\cup_{I\in B}I\ne\emp$; let $j$ be in this last set.
For some $s\ge0$ we have $\k^s(j)\ne\aa$, so that if $B'=\k^s(B)$ we have
$\aa\n\ov{[1,D]}-\cup_{I\in B'}I$, that is $\cup_{I\in B'}I\sub[1,D]$ so that $I^*=I$ for any
$I\in B'$. It follows that $B'{}^*=B'$. By (a) we have $B'\in\SS^*_D$ hence $B'\in\SS_D$. Since
$B'\in\SS_D$, $B'=B'{}^1$ we have by definition $B'\in S_D$. This proves (b).

\mpb

In view of 1.11, (b) provides a non-inductive description of $\SS^*_D$ hence also of $\SS_D$,
which is simpler that that in \cite{L20, 1.3(c)}.

\head 4. Tables in classical types\endhead
\subhead 4.1\endsubhead
Given $\tB\in\SS^*_D$ we can define $B\in S_D$ by $\l\i(\tB)=(B,k)$ (see 2.8) and then
$E=\EE_B\in\cf(V_D)$ and $X=(E^0)^!\in\cc(V_D^1)$, $Y=E^1\in\cc(V_D^1)$; we can also form
$\e_D(\tB)\in V_D$.
We record the assignments $\tB\m(Y,X)$ and $\tB\m\e_D(\tB)$ in the tables of
4.2-4.4 (representing the cases where $D=2,4,6$) and in the tables of 4.5-4.8
(representing the cases where $D=1,3,5,7$).

The table for $V_D$ consists of several subtables, one for each $X\in\cc(V_D^1)$.
The subtable indexed by $X$ has the name $X$ in a box, has one row for each $Y\in\cc(V_D^1)$
such that $(Y,X)\in\ovm{\cc}(V_D^1)$ and has a list of the elements $\tB\in\SS^*_D$ which give
rise as above to $Y\sub X$; the image $\e_D(\tB)\in V_D$ of such a $\tB$ is also given.

Given $\tB\in\SS'{}^*_D$ with $D$ odd we can define $(Y,X)$ as above; we set
$Y'=\p(Y)\in\cc(V'_D{}^1),X'=\p(X)\in\cc(V'_D{}^1)$.
We record the assignments $\tB\m(Y',X')$ and $\tB\m\e'_D(\tB)$ in the tables of 4.9-4.12
(representing the cases where $D=1,3,5,7$).
The table for $V'_D$ consists of several subtables, one for each $X'\in\cc(V'_D{}^1)$. The
subtable indexed by $X'$ has the name $X'$ in a box, has one row for each $Y'\in\cc(V'_D{}^1)$
such that $(Y',X')\in\ovm{\cc}(V'_D{}^1)$ and has a list of the elements $\tB\in\SS'{}^*_D$ which
give rise as above to $Y'\sub X'$; the image $\e'_D(\tB)\in V'_D$ of such a $\tB$ is also given.

We now explain the notation in the tables.
Any $\tB$ in $\SS^*_D$ (or in $\SS'{}^*_D$ with $D$ odd) is represented as list of intervals.
For example $(6,56\aa,456\aa)$ (for $V_6$) represents the set of intervals
$\{6\},\{5,6,\aa\}, \{4,5,6,\aa\}$. The elements of $V_D$ are written as $i_1i_2\do i_k$
instead of $e_{i_1}+e_{i_2}+\do+e_{i_k}$. For example $236$ represents $e_2+e_3+e_6\in V_6$.
The elements of $V'_D$ are written as $i_1i_2\do i_k$ instead of
$\p(e_{i_1}+e_{i_2}+\do+e_{i_k})$. 
The subspaces $X$ or $Y$ of $V_D^1$ are represented by a sequence of generating vectors.
For example $(13,5)$ represents the subspace spanned by $e_1+e_3$ and $e_5$.
The subspaces $X'$ or $Y'$ of $V'_D{}^1$ are written in the form $\p(X)$ or $\p(Y)$
with $X$ or $Y$ written as above.

\subhead 4.2. Table for $V_2$\endsubhead

$$\boxed{X=(1)}$$   

$\emp\m0$; $(\aa)\m12$; $Y=0$

$(1)\m1$;            $Y=(1)$

\mpb

$$\boxed{X=0}$$

$(2)\m2$;  $Y=0$

\subhead 4.3. Table for $V_4$\endsubhead

$$\boxed{X=(1,3)}$$

$\emp\m0$; $(\aa)\m1234$; $(\aa,1\aa4)\m23$;   $Y=0$

$(1)\m1$; $(1,\aa12)\m34$;                 $Y=(1)$

$(3)\m3$; $(3,\aa)\m124$;                  $Y=(3)$

$(1,3)\m13$;                            $Y=(1,3)$

\mpb

$$\boxed{X=(13)}$$

$(2)\m2$; $(2,\aa)\m134$;                  $Y=0$

$(2,123)\m123$;                          $Y=(13)$

\mpb

$$\boxed{X=(1)}$$

$(4)\m4$; $(4,\aa43)\m12$;                  $Y=0$

$(1,4)\m14$;                               $Y=(1)$

\mpb

$$\boxed{X=(3)}$$

$(3,234)\m234$;                            $Y=(3)$

\mpb

$$\boxed{X=0}$$

$(2,4)\m24$;                              $Y=0$

\subhead 4.4. Table for $V_6$\endsubhead

$$\boxed{X=(1,3,5)}$$  

$\emp\m0$; $(\aa)\m123456$; $(\aa,1\aa6)\m2345$; $(\aa,1\aa6,21\aa65)\m1256$;   $Y=0$

$(1)\m1$; $(1,\aa12)\m3456$; $(1,21\aa,321\aa6)\m145$;                    $Y=(1)$

$(3)\m3$; $(3,\aa)\m12456$; $(3,\aa,1\aa6)\m245$;                         $Y=(3)$

$(5)\m5$; $(5,\aa)\m12346$; $(5,\aa,1\aa654)\m23$;                        $Y=(5)$

$(1,3)\m13$; $(1,3,\aa1234)\m56$;                                     $Y=(1,3)$

$(1,5)\m15$; $(1,5,\aa12)\m346$;                                      $Y=(1,5)$

$(3,5)\m35$; $(3,5,\aa)\m1246$;                                       $Y=(3,5)$

$(1,3,5)\m135$;                                                     $Y=(1,3,5)$

\mpb

$$\boxed{X=(13,5)}$$

$(2)\m2$; $(2,\aa)\m13456$; $(2,\aa,321\aa6)\m45$;                        $Y=0$

$(2,5)\m25$; $(2,5,\aa)\m1346$;                                       $Y=(5)$

$(2,123)\m123$; $(2,123,\aa1234)\m256$;                               $Y=(13)$

$(2,5,123)\m1235$;                                                  $Y=(13,5)$

\mpb

$$\boxed{X=(1,35)}$$

$(4)\m4$; $(4,\aa)\m12356$; $(4,\aa,1\aa6)\m235$;                        $Y=0$

$(1,4)\m14$; $(1,4,\aa12)\m356$;                                     $Y=(1)$

$(4,345)\m345$; $(4,345,\aa)\m126$;                                  $Y=(35)$

$(1,4,345)\m1345$;                                                 $Y=(1,35)$

\mpb

$$\boxed{X=(1,3)}$$

$(6)\m6$; $(6,56\aa)\m1234$; $(6,56\aa, 1\aa654)\m236$;                    $Y=0$

$(1,6)\m16$; $(1,6,21\aa65)\m34$;                                    $Y=(1)$

$(3,6)\m36$; $(3,6,56\aa)\m124$;                                     $Y=(3)$

$(1,3,6)\m136$;                                                    $Y=(1,3)$

\mpb

$$\boxed{X=(1,5)}$$

$(5,456)\m456$; $(5,456,3456\aa)\m125$;                              $Y=(5)$

$(1,5,456)\m1456$;                                                 $Y=(1,5)$

\mpb

$$\boxed{X=(3,5)}$$

$(3,5,23456)\m23456$;                                                    $Y=(3,5)$

\mpb

$$\boxed{X=(3,135)}$$

$(3,234)\m234$; $(3,234,\aa)\m156$;                                      $Y=(3)$

$(3,234,12345)\m1245$;                                                 $Y=(3,135)$

\mpb

$$\boxed{X=(135)}$$

$(2,4)\m24$; $(2,4,\aa)\m1356$;                                          $Y=0$

$(2,4,12345)\m12345$;                                                  $Y=(135)$

\mpb

$$\boxed{X=(1)}$$

$(4,6)\m46$; $(4,6,3456\aa)\m12$;                                        $Y=0$

$(1,4,6)\m146$;                                                        $Y=(1)$

\mpb

$$\boxed{X=(3)}$$

$(3,6,234)\m2346$;                                                     $Y=(3)$

\mpb

$$\boxed{X=(5)}$$

$(2,5,456)\m2456$;                                                     $Y=(5)$

\mpb

$$\boxed{X=(13)}$$

$(2,6)\m26$; $(2,6,56\aa)\m134$;                                         $Y=0$

$(2,6,123)\m1236$;                                                     $Y=(13)$

\mpb

$$\boxed{X=(35)}$$

$(4,345,23456)\m2356$;                                                 $Y=(35)$

\mpb

$$\boxed{X=0}$$

$(2,4,6)\m246$;                                                        $Y=0$

\subhead 4.5. Table for $V_1$\endsubhead

$$\boxed{X=(1)}$$

$\emp\m0$;    $Y=0$

$(1)\m1$;   $Y=(1)$

\subhead 4.6. Table for $V_3$\endsubhead

$$\boxed{X=(1,3)}$$

$\emp\m0$; $(\aa3)\m12$; $(1\aa)\m23$;  $Y=0$   

$(1)\m1$;      $Y=(1)$

$(3)\m3$;      $Y=(3)$

$(1,3)\m13$;   $Y=(1,3)$

\mpb

$$\boxed{X=(13)}$$

$(2)\m2$;  $Y=0$

$(2,123)\m123$;  $Y=(13)$

\subhead 4.7. Table for $V_5$\endsubhead

$$\boxed{X=(1,3,5)}$$

$\emp\m0$; $(\aa5)\m1234$; $(1\aa)\m2345$;     $Y=0$

$(1)\m1$; $(1,21\aa5)\m34$; $(1,321\aa)\m145$;   $Y=(1)$

$(3)\m3$; $(3,\aa5)\m124$; $(3,1\aa)\m245$;      $Y=(3)$

$(5)\m5$; $(5,1\aa54)\m23$; $(5,\aa543)\m125$;   $Y=(5)$

$(1,3)\m13$;                                 $Y=(1,3)$

$(1,5)\m15$;                                 $Y=(1,5)$

$(3,5)\m35$;                                 $Y=(3,5)$

$(1,3,5)\m135$;                              $Y=(1,3,5)$

\mpb

$$\boxed{X=(13,5)}$$

$(2)\m2$; $(2,\aa5)\m134$; $(2,321\aa)\m45$;  $Y=0$

$(2,5)\m25$;                              $Y=(5)$        

$(2,123)\m123$;                           $Y=(13)$

$(2,5,123)\m1235$;                        $Y=(13,5)$

\mpb

$$\boxed{X=(1,35)}$$

$(4)\m4$;$(4,\aa543)\m12$; $(4,1\aa)\m235$;   $Y=0$

$(1,4)\m14$;                              $Y=(1)$

$(4,345)\m345$;                           $Y=(35)$

$(1,4,345)\m1345$;                        $Y=(1,35)$

\mpb

$$\boxed{X=(3,135)}$$

$(3,234)\m234$;                     $Y=(3)$

$(3,234,12345)\m1245$;              $Y=(3,135)$

\mpb

$$\boxed{X=(135)}$$

$(2,4)\m24$;                       $Y=0$

$(2,4,12345)\m12345$;              $Y=(135)$

\subhead 4.8. Table for $V_7$\endsubhead

$$\boxed{X=(1,3,5,7)}$$

$\emp\m0$; $(\aa7)\m123456$; $(1\aa)\m234567$; $(\aa7,1\aa76,21\aa765)\m1256$;

                       $(1\aa,21\aa7,321\aa76)\m2367$;   $Y=0$

$(1)\m1$; $(1,21\aa7)\m3456$; $(1,321\aa)\m14567$;           $Y=(1)$

$(3)\m3$; $(3,\aa7)\m12456$; $(3,1\aa)\m24567$;              $Y=(3)$

$(5)\m5$; $(5,\aa7)\m12346$; $(5,1\aa)\m23467$;              $Y=(5)$

$(7)\m7$; $(7,1\aa76)\m2345$; $(7,\aa765)\m12347$;           $Y=(7)$

$(1,3)\m13$; $(1,3,4321\aa7)\m56$; $(1,3,54321\aa)\m1367$;   $Y=(1,3)$

$(1,5)\m15$; $(1,5,21\aa7)\m346$; $(1,5,321\aa)\m1467$;      $Y=(1,5)$

$(1,7)\m17$; $(1,7,321\aa76)\m145$; $(1,7,21\aa765)\m347$;   $Y=(1,7)$

$(3,5)\m35$; $(3,5,\aa7)\m1246$; $(3,5,1\aa)\m2467$;         $Y=(3,5)$

$(3,7)\m37$; $(3,7,1\aa76)\m245$; $(3,7,\aa765)\m1247$;      $Y=(3,7)$

$(5,7)\m57$; $(5,7,1\aa7654)\m23$; $(5,7,\aa76543)\m1257$;   $Y=(5,7)$

$(1,3,5)\m135$;                                          $Y=(1,3,5)$

$(1,3,7)\m137$;                                          $Y=(1,3,7)$

$(1,5,7)\m157$;                                          $Y=(1,5,7)$

$(3,5,7)\m357$;                                          $Y=(3,5,7)$

$(1,3,5,7)\m1357$;                                       $Y=(1,3,5,7)$

\mpb

$$\boxed{X=(13,5,7)}$$

$(2)\m2$; $(2,\aa7)\m13456$; $(2,321\aa)\m4567$;       $Y=0$

$(2,5)\m25$; $(2,5,\aa7)\m1346$; $(2,5,321\aa)\m467$;   $Y=(5)$

$(2,7)\m27$; $(2,7,321\aa76)\m45$; $(2,7,\aa765)\m1347$;   $Y=(7)$

$(2,123)\m123$; $(2,123,4321\aa7)\m256$; $(2,123,54321\aa)\m12367$;  $Y=(13)$

$(2,5,7)\m257$;                                                  $Y=(5,7)$

$(2,5,123)\m1235$;                                               $Y=(13,5)$

$(2,7,123)\m1237$;                                               $Y=(13,7)$

$(2,5,7,123)\m12357$;                                            $Y=(13,5,7)$

\mpb

$$\boxed{X=(1,35,7)}$$

$(4)\m4$; $(4,\aa7)\m12356$; $(4,1\aa)\m23567$;    $Y=0$

$(1,4)\m14$; $(1,4,21\aa7)\m356$; $(1,4,54321\aa)\m167$; $Y=(1)$

$(4,7)\m47$; $(4,7,1\aa76)\m235$; $(4,7,\aa76543)\m127$; $Y=(7)$

$(4,345)\m345$; $(4,345,\aa7)\m126$; $(4,345,1\aa)\m267$; $Y=(35)$

$(1,4,7)\m147$;                                       $Y=(1,7)$

$(1,4,345)\m1345$;                                    $Y=(1,35)$

$(4,7,345)\m3457$;                                    $Y=(35,7)$

$(1,4,7,345)\m13457$;                                 $Y=(1,35,7)$

\mpb

$$\boxed{X=(1,3,57)}$$

$(6)\m6$; $(6,\aa765)\m1234$; $(6,1\aa)\m23457$;       $Y=0$

$(1,6)\m16$; $(1,6,21\aa765)\m34$; $(1,6,321\aa)\m1457$;     $Y=(1)$

$(3,6)\m36$; $(3,6,\aa765)\m124$; $(3,6,1\aa)\m2457$;        $Y=(3)$

$(6,567)\m567$; $(6,567,1\aa7654)\m236$; $(6,567,\aa76543)\m12567$; $Y=(57)$

$(1,3,6)\m136$;                                                 $Y=(1,3)$

$(1,6,567)\m1567$;                                              $Y=(1,57)$

$(3,6,567)\m3567$;                                              $Y=(3,57)$

$(1,3,6,567)\m13567$;                                           $Y=(1,3,57)$

\mpb

$$\boxed{X=(3,135,7)}$$

$(3,234)\m234$; $(3,234,\aa7)\m156$; $(3,234,54321\aa)\m367$;  $Y=(3)$

$(3,7,234)\m2347$;                                         $Y=(3,7)$

$(3,234,12345)\m1245$;                                     $Y=(3,135)$

$(3,7,234,12345)\m12457$;                                  $Y=(3,135,7)$

\mpb

$$\boxed{X=(1,5,357)}$$

$(5,456)\m456$; $(5,456,\aa76543)\m125$; $(5,456,1\aa)\m237$;  $Y=(5)$

$(1,5,456)\m1456$;                                         $Y=(1,5)$

$(5,456,34567)\m3467$;                                     $Y=(5,357)$

$(1,5,456,34567)\m13467$;                                   $Y=(1,5,357)$

\mpb

$$\boxed{X=(3,5,1357)}$$

$(3,5,23456)\m23456$;                                              $Y=(3,5)$

$(3,5,23456,1234567)\m12467$;                                      $Y=(3,5,1357)$

\mpb

$$\boxed{X=(135,7)}$$

$(2,4)\m24$; $(2,4,\aa7)\m1356$; $(2,4,54321\aa)\m67$;    $Y=0$

$(2,4,7)\m247$;                                       $Y=(7)$

$(2,4,12345)\m12345$;                                 $Y=(135)$

$(2,4,7,12345)\m123457$;                              $Y=(135,7)$

\mpb

$$\boxed{X=(13,57)}$$

$(2,6)\m26$; $(2,6,\aa765)\m134$; $(2,6,321\aa)\m457$;    $Y=0$

$(2,6,123)\m1236$;                                    $Y=(13)$

$(2,6,567)\m2567$;                                    $Y=(57)$

$(2,6,123,567)\m123567$;                              $Y=(13,57)$

\mpb

$$\boxed{X=(1,357)}$$

$(4,6)\m46$; $(4,6,\aa76543)\m12$; $(4,6,1\aa)\m2357$;    $Y=0$

$(1,4,6)\m146$;                                       $Y=(1)$

$(4,6,34567)\m34567$;                                 $Y=(357)$

$(1,4,6,34567)\m134567$;                              $Y=(1,357)$

\mpb

$$\boxed{X=(3,1357)}$$

$(3,6,234)\m2346$;                                     $Y=(3)$

$(3,6,234,1234567)\m124567$;                          $Y=(3,1357)$

\mpb

$$\boxed{X=(5,1537)}$$

$(2,5,456)\m2456$;                                       $Y=(5)$

$(2,5,456,1234567)\m123467$;                             $Y=(5,1357)$

\mpb

$$\boxed{X=(35,1357)}$$

$(4,345,23456)\m2356$;                                            $Y=(35)$

$(4,345,23456,1234567)\m1267$;                                    $Y=(35,1357)$

\mpb

$$\boxed{X=(1357)}$$

$(2,4,6)\m246$;                                                   $Y=0$

$(2,4,6,1234567)\m1234567$;                                       $Y=(1357)$

\subhead 4.9. Table for $V'_1$\endsubhead

$$\boxed{X'=\p(1)}$$

$\emp\m0$;    $Y'=0$

\subhead 4.10. Table for $V'_3$\endsubhead

$$\boxed{X'=\p(1,3)}$$

$\emp\m0$; $(\aa3)\m12$; $Y'=0$

$(1)\m1$;  $Y'=\p(1)$

\mpb

$$\boxed{X'=\p(13)}$$

$(2)\m2$;   $Y'=0$

\subhead 4.11. Table for $V'_5$\endsubhead

$$\boxed{X'=\p(1,3,5)}$$

$\emp\m0$; $(\aa5)\m1234$;  $Y'=0$

$(1)\m1$; $(1,21\aa5)\m34$;  $Y'=\p(1)$

$(3)\m3$; $(3,\aa5)\m124$;   $Y'=\p(3)$

$(5,1\aa54)\m23$;        $Y'=\p(5)$

$(1,3)\m13$;               $Y'=\p(1,3)$

\mpb

$$\boxed{X'=\p(13,5)}$$

$(2)\m2$; $(2,\aa5)\m134$;  $Y'=0$

$(2,123)\m123$;           $Y'=\p(13)$

\mpb

$$\boxed{X'=\p(1,35)}$$

$(4)\m4$; $(4,\aa543)\m12$;  $Y'=0$

$(1,4)\m14$;               $Y'=\p(1)$

\mpb

$$\boxed{X'=\p(3,135)}$$

$(3,234)\m234$;            $Y'=\p(3)$

\mpb

$$\boxed{X'=\p(135)}$$

$(2,4)\m24$;               $Y'=0$

\subhead 4.12. Table for $V'_7$\endsubhead

$$\boxed{X'=\p(1,3,5,7)}$$

$\emp\m0$; $(\aa7)\m123456$; $(\aa7,1\aa76,21\aa765)\m1256$; $Y'=0$

$(1)\m1$; $(1,21\aa7)\m3456$;                          $Y'=\p(1)$

$(3)\m3$; $(3,\aa7)\m12456$;                           $Y'=\p(3)$

$(5)\m5$; $(5,\aa7)\m12346$;                           $Y'=\p(5)$

$(7,1\aa76)\m2345$;                         $Y'=\p(7)$

$(1,3)\m13$; $(1,3,4321\aa7)\m56$;                     $Y'=\p(1,3)$

$(1,5)\m15$; $(1,5,21\aa7)\m346$;                      $Y'=\p(1,5)$

$(1,7,321\aa76)\m145$;                         $Y'=\p(1,7)$

$(3,5)\m35$; $(3,5,\aa7)\m1246$;                       $Y'=\p(3,5)$

$(3,7,1\aa76)\m245$;                         $Y'=\p(3,7)$

$(5,7,1\aa7654)\m23$;                         $Y'=\p(5,7)$

$(1,3,5)\m135$;                                      $Y'=\p(1,3,5)$

\mpb

$$\boxed{X'=\p(13,5,7)}$$

$(2)\m2$; $(2,\aa7)\m13456$;      $Y'=0$

$(2,5)\m25$; $(2,5,\aa7)\m1346$;  $Y'=\p(5)$

$(2,7,321\aa76)\m45$;   $Y'=\p(7)$

$(2,123)\m123$; $(2,123,4321\aa7)\m256$;  $Y'=\p(13)$

$(2,5,123)\m1235$;                      $Y'=\p(13,5)$

\mpb

$$\boxed{X'=\p(1,35,7)}$$

$(4)\m4$; $(4,\aa7)\m12356$;  $Y'=0$

$(1,4)\m14$; $(1,4,21\aa7)\m356$;  $Y'=\p(1)$

$(4,7,1\aa76)\m235$;    $Y'=\p(7)$

$(4,345)\m345$; $(4,345,\aa7)\m126$;   $Y'=\p(35)$

$(1,4,345)\m1345$;                   $Y'=\p(1,35)$

\mpb

$$\boxed{X'=\p(1,3,57)}$$

$(6)\m6$; $(6,\aa765)\m1234$;      $Y'=0$

$(1,6)\m16$; $(1,6,21\aa765)\m34$;  $Y'=\p(1)$

$(3,6)\m36$; $(3,6,\aa765)\m124$;   $Y'=\p(3)$

$(6,567,1\aa7654)\m236$;      $Y'=\p(57)$   

$(1,3,6)\m136$;                   $Y'=\p(1,3)$

\mpb

$$\boxed{X'=\p(3,135,7)}$$

$(3,234)\m234$; $(3,234,\aa7)\m156$;   $Y'=\p(3)$

$(3,234,12345)\m1245$;               $Y'=\p(3,135)$

\mpb

$$\boxed{X'=\p(1,5,357)}$$

$(5,456)\m456$; $(5,456,\aa76543)\m125$;   $Y'=\p(5)$

$(1,5,456)\m1456$;                       $Y'=\p(1,5)$

\mpb

$$\boxed{X'=\p(3,5,1357)}$$

$(3,5,23456)\m23456$;   $Y'=\p(3,5)$

\mpb

$$\boxed{X'=\p(135,7)}$$

$(2,4)\m24$; $(2,4,\aa7)\m1356$;      $Y'=0$

$(2,4,12345)\m12345$;               $Y'=\p(135)$

\mpb

$$\boxed{X'=\p(13,57)}$$

$(2,6)\m26$; $(2,6,\aa765)\m134$;    $Y'=0$

$(2,6,123)\m1236$;                $Y'=\p(13)$

\mpb

$$\boxed{X'=\p(1,357)}$$

$(4,6)\m46$; $(4,6,\aa76543)\m12$;    $Y'=0$

$(1,4,6)\m146$;                     $Y'=\p(1)$

\mpb

$$\boxed{X'=\p(3,1357)}$$

$(3,6,234)\m2346$;                 $Y'=\p(3)$

\mpb

$$\boxed{X'=\p(5,1537)}$$

$(2,5,456)\m2456$;                $Y'=\p(5)$

\mpb

$$\boxed{X'=\p(35,1357)}$$

$(4,345,23456)\m2356$;                                               $Y'=\p(35)$

\mpb

$$\boxed{X'=\p(1357)}$$

$(2,4,6)\m246$;                  $Y'=0$

\head 5. Exceptional types\endhead
\subhead 5.1\endsubhead
We now return to the setup in 0.1 and assume that $G$ is of exceptional type.
Recall that we have fixed a family $c$ in $\Irr(W)$.
We must be in one of the following cases.

(i) $|c|=1$, $\cg_c=S_1$.

(ii) $|c|=2$ (with $W$ of type $E_7$ or $E_8$), $\cg_c=S_2$.

(iii) $|c|=3$, $\cg_c=S_2$.

(iv) $|c|=4$ (with $W$ of type $G_2$), $\cg_c=S_3$.

(v) $|c|=5$ (with $W$ of type $E_6,E_7$ or $E_8$), $\cg_c=S_3$.

(vi) $|c|=11$ (with $W$ of type $F_4$), $\cg_c=S_4$.

(vii) $|c|=17$ (with $W$ of type $E_8$), $\cg_c=S_5$.

Here for $n\in[1,5]$, $S_n$ denotes the group of permutations of $[1,n]$.

\subhead 5.2\endsubhead
In this subsection we assume that $\cg_c=S_5$. We write $H_{5!}=S_5$.
Let $H_{4!}$ be the group of all $\s\in S_5$ which
map $[1,4]$ to itself and $5$ to itself.
Let $H_{2!3!}$ be the group of all $\s\in S_5$ which map $[1,2]$
to itself and $[3,5]$ to itself. Let $H_8$ be the
the group of all $\s\in S_5$ which commute with the permutation
$1\m2\m1,3\m4\m3$, $5\m5$.
Let $H_{3!}$ be the group of all $\s\in S_5$ which map $[3,5]$ to
itself and $1$ to $1$, $2$ to $2$.
Let $H_{2!2!}$ be the group of all $\s\in S_5$ which map $[1,2]$
to itself, $[3,4]$ to itself and $5$ to $5$.
Let $H_{2!}$ be the subgroup of all $\s\in S_5$ which map
$[1,2]$ to $[1,2]$, $3$ to $3$, $4$ to $4$ and $5$ to $5$.
Let $H_{1!}=\{1\}\sub S_5$. Let
$$\HH_c=\{H_{5!},H_{4!},H_{2!3!},H_8,H_{2!2!},H_{3!},H_{2!},H_{1!}\}.$$
We define $\ovm{\HH}_c$ as the set of all pairs $(H_i,H_j)$
where $i\le j$ are such that $H_i$ is a normal subgroup of $H_j$ but $(i,j)\ne(1,8)$.

\subhead 5.3\endsubhead
In this subsection we assume that $\cg_c=S_4$. We write $H_{4!}=S_4$.
Let $H_8$ be the the group of all $\s\in S_4$ which commute with the permutation
$1\m2\m1,3\m4\m3$. Let $H_{3!}$ be the group of all $\s\in S_4$ which map $[1,3]$ to
itself and $4$ to $4$.
Let $H_{2!2!}$ be the group of all $\s\in S_4$ which map $[1,2]$
to itself, $[3,4]$ to itself.
Let $H_{2!}$ be the subgroup of all $\s\in S_4$ which map
$[1,2]$ to $[1,2]$, $3$ to $3$, $4$ to $4$.
Let $H_{1!}=\{1\}\sub S_4$. Let
$$\HH_c=\{H_{4!},H_8,H_{2!2!},H_{3!},H_{2!},H_{1!}.\}$$
We define $\ovm{\HH}_c$ as the set of all pairs $(H_i,H_j)$
where $i\le j$ are such that $H_i$ is a normal subgroup of $H_j$ but $(i,j)\ne(1,8)$.

\subhead 5.4\endsubhead
In this subsection we assume that $\cg_c=S_3$. We write $H_{3!}=S_3$.
Let $H_{2!}$ be the subgroup of all $\s\in S_3$ which map
$[1,2]$ to $[1,2]$, $3$ to $3$. Let $H_{1!}=\{1\}\sub S_3$. Let
$$\HH_c=\{H_{3!},H_{2!},H_{1!}.\}$$
We define $\ovm{\HH}_c$ as the set of all pairs $(H_i,H_j)$
where $i\le j$ are such that $H_i$ is a normal subgroup of $H_j$.

\subhead 5.5\endsubhead
If $\cg_c=S_2$ let $H_{2!}=S_2,H_{1!}=\{1\}\sub S_2$. Let $\HH_c=\{H_{2!},H_{1!}\}$.
We define $\ovm{\HH}_c$ as the set of all pairs $H_i\sub H_j$ where $i\le j$.
If $\cg_c=S_1$ let $H_{1!}=\{1\}=S_1$. Let $\HH_c=\{H_{1!}\}$.
We define $\ovm{\HH}_c$ as the set consisting of $(H_{1!},H_{1!})$.

\subhead 5.6\endsubhead
In each of the cases in 5.2-5.5, for any $(H_i,H_j)$ in $\ovm{\HH}_c$ we can describe
$H_j/H_i$ as follows.

If $i=j$ then $H_j/H_i=\{1\}$.

If $i=k!$ then $H_i/\{1\}=S_k$ canonically.

If $i=2!3!$ then $H_i/H_{2!}=S_3,H_i/H_{3!}=S_2,H_i/\{1\}=S_2\T S_3$ canonically.

If $i=2!2!$ then $H_i/H_{2!}=S_2,H_i/\{1\}=S_2\T S_2$ canonically.

If $i=8$ then $H_i/H_{2!2!}=S_2$ canonically.

\subhead 5.7\endsubhead
When $n>1$ let $\mu'_n$ be the subset of $\CC^*$ consisting of the primitive $n$-th roots of $1$.

Let $\G$ be one of the groups $S_n,n=1,2,3,4,5$ or $S_2\T S_2$ or $S_3\T S_2$.
We define a subset $Prim(\G)$ of $\CC[M(\G)]$ as follows.
If $\G=S_1$ then $Prim(\G)$ consists of $(1,1)$.

If $\G=S_n$, $n=2,3,4,5$, then $Prim(\G)$ consists of $(1,1)$ and of $\L_e$, $e\in\mu'_n$
(as in \cite{L20, 3.2}).

If $\G=S_2\T S_3$, then $Prim(\G)$ consists of the five elements $\L_{-1}\bxt(1,1)$,
$\L_{-1}\bxt\L_e,(1,1)\bxt\L_e$ ($e\in\mu'_3$) of $\CC[M(S_2)]\ot\CC[ M(S_3)]=\CC[M(\G)]$ and of
$(1,1)\in\CC[M(\G)]$. (Note that in this case we have
$Prim(\G)=Prim(S_2)\ot Prim(S_3)\sub\CC[M(S_2)]\ot\CC[M(S_3)]=\CC[M(\G)]$.)

If $\G=S_2\T S_2$, then $Prim(\G)$ consists of the two elements
$\L_{-1}\bxt\L_{-1},\L_{-1}\bxt(1,1)$ of $\CC[M(S_2)]\ot\CC[M(S_2)]=\CC[M(\G)]$ and of
$(1,1)\in\CC[M(\G)]$.
(Note that the two factors in $S_2\T S_2$ play an asymmetric role: in fact, $S_2\T S_2$ can be
viewed as a two-dimensional vector space with a given ordered basis.)

\subhead 5.8\endsubhead
In each of the cases in 5.2-5.5, for any $(H_i,H_j)$ in $\ovm{\HH}_c$ we set

$Prim(H_i,H_j)=Prim(H_j/H_i)$;
\nl
we use the identifications in 5.6 and the definitions in 5.7. 
Therefore both the source and target of $\Th$ in 0.6(a) are defined. The existence and uniqueness
of such a $\Th$ has been already verified in \cite{L19}, \cite{L20}. (Note that in
\cite{L20} the definition of $Prim(S_5)$ and that of $Prim(S_3)$ for $G$ of type $G_2$
is different from the present one, but this does not affect the proofs.) Thus the statements in
0.6 and hence those in 0.2 are established when $G$ is of exceptional type.

In 5.9-5.13 we describe the bijection 0.6(a) in each of the cases in 5.2-5.5. We give a table
for each type of $\cg_c$. Each table consists of several subtables, one for each $H'\in\HH_c$.
The subtable indexed by $H'\in\HH_c$ has the name $X=H'$ (in a box), has one row for each
$Y=H\in\HH_c$ such that $(H,H')\in\ovm{\HH}_c$ and has a list of the elements of
$M(\cg_c)$ (in the notation of \cite{L84, \S4}) which are in the fibre of $\a_c$ (in 0.3) at $(H,H')$.

\subhead 5.9. Table for $\cg_c=S_1$ \endsubhead

$$\boxed{X=H_{1!}}$$

$(1,1)$;  $Y=1$

\subhead 5.10. Table for $\cg_c=S_2$ \endsubhead

$$\boxed{X=H_{2!}}$$

$(1,1)$,    $(g_2,\e)$;   $Y=1$

$(g_2,1)$;   $Y=H_{2!}$    

\mpb

$$\boxed{X=1}$$

$(1,\e)$;   $Y=1$

\subhead 5.11. Table for $\cg_c=S_3$ \endsubhead

$$\boxed{X=H_{3!}}$$

$(1,1)$, $(g_3,\th),(g_3,\th^2)$;   $Y=1$

$(g_3,1)$;   $Y=H_{3!}$

\mpb

$$\boxed{X=H_{2!}}$$

$(1,r)$,    $(g_2,\e)$;   $Y=1$

$(g_2,1)$;   $Y=H_{2!}$    

\mpb

$$\boxed{X=1}$$

$(1,\e)$;   $Y=1$

\subhead 5.12. Table for $\cg_c=S_4$ \endsubhead

$$\boxed{X=H_{4!}}$$

$(1,1)$; $(g_4,i),g_4(-i)$;                          $Y=1$

$(g_4,1)$;                                     $Y=H_{4!}$

\mpb

$$\boxed{X=H_8}$$

$(g'_2,1),(g_4,-1)$                             $Y=H_{2!2!}$

$(g'_2,\e')$;                             $Y=H_8$

\mpb

$$\boxed{X=H_{2!2!}}$$

$(1,\s)$, $(g_2,\e'),   (g'_2,\e)$; $Y=1$

$(g_2,1)$,  $(g'_2,r)$;   $Y=H_{2!}$

$(g'_2,\e'')$;           $Y=H_{2!2!}$

\mpb

$$\boxed{X=H_{3!}}$$

$(1,\l^1)$, $(g_3,\th),(g_3,\th^2)$;   $Y=1$

$(g_3,1)$;   $Y=H_{3!}$

\mpb

$$\boxed{X=H_{2!}}$$

$(1,\l^2)$,    $(g_2,\e)$;   $Y=1$

$(g_2,\e'')$;   $Y=H_{2!}$    
                                                                                 
\mpb

$$\boxed{X=1}$$

$(1,\l^3)$;   $Y=1$

\subhead 5.13. Table for $\cg_c=S_5$ \endsubhead
$$\boxed{X=H_{5!}}$$
$(1,1)$; $(g_5,\z),(g_5,\z^2),(g_5,\z^3),(g_5,\z^4)$;  $Y=1$

$(g_5,1)$;                            $Y=H_{5!}$

\mpb

$$\boxed{X=H_{4!}}$$

$(1,\l^1)$; $(g_4,i),g_4(-i)$;                          $Y=1$

$(g_4,1)$;                                     $Y=H_{4!}$

\mpb

$$\boxed{X=H_{2!3!}}$$

$(1,\nu)$, $(g_2,-1),(g_3,\th),(g_3,\th^2),(g_6,-\th),(g_6,-\th^2)$; $Y=1$
 
$(g_2,1),(g_6,\th),(g_6,\th^2)$; $Y=H_{2!}$

$(g_3,1),(g_6,-1)$; $Y=H_{3!}$

$(g_6,1)$; $Y=H_{2!3!}$

\mpb

$$\boxed{X=H_8}$$

$(g'_2,1),(g_4,-1)$                             $Y=H_{2!2!}$

$(g'_2,\e')$;                             $Y=H_8$

\mpb

$$\boxed{X=H_{2!2!}}$$

$(1,\nu')$, $(g_2,-r),   (g'_2,\e)$; $Y=1$

$(g_2,r)$,  $(g'_2,r)$;   $Y=H_{2!}$

$(g'_2,\e'')$;           $Y=H_{2!2!}$

\mpb

$$\boxed{X=H_{3!}}$$

$(1,\l^2)$, $(g_3,\e\th),(g_3,\e\th^2)$;   $Y=1$

$(g_3,\e)$;   $Y=H_{3!}$

\mpb

$$\boxed{X=H_{2!}}$$

$(1,\l^3)$,    $(g_2,-\e)$;   $Y=1$

$(g_2,\e)$;   $Y=H_{2!}$

\mpb

$$\boxed{X=1}$$

$(1,\l^4)$;   $Y=1$

\head 6. The partition $\cu_c=\sqc_{\gg\in\un\PP}\cu_c^\gg$\endhead
\subhead 6.1\endsubhead
Assume that $G$ is as in 0.6(i).
The partition 0.2(a) of $\cu_c$ (or equivalently of $M(\cg_c)=V_D$, $D$ even)
corresponds under the bijection $\e_D$ in 3.3(a) to the partition
$$\SS^*_D=\sqc_{k\in\NN}\SS^*_D(k)$$
where $\SS^*_D(k)$ is as in 3.9.
From the proof in 3.7 we see that 0.6(b) and 0.7(a) hold in our case and that
if $\gg\in\un\PP$ corresponds to $k$ as above, then

(a) $\ovm{\HH}_c^\gg=\{(Y\sub X)\in\ovm{\cc}(V_D^1);k\le\dim(X/Y)\}$.

\subhead 6.2\endsubhead
Assume that $G$ is as in 0.6(ii).
The partition 0.2(a) of $\cu_c$ (or equivalently of $M(\cg_c)=V'_D$, $D$ odd)
corresponds under the bijection $\e'_D$ in 3.3(b) to the partition
$$\SS'{}^*_D=\sqc_{k\in\{0,1,3,5,\do\}}\SS'{}^*_D(k)$$
where $\SS'{}^*_D(k)$ is as in 3.10.
From the proof in 3.8 we see that 0.6(b) and 0.7(a) hold in our case and that
if $\gg\in\un\PP$ corresponds to $k$ as above then

(a) $\ovm{\HH}_c^\gg=\{(Y'\sub X')\in\ovm{\cc}(V'_D{}^1);k\le\dim(X'/Y')\}$ if $k>0$;

(b) $\ovm{\HH}_c^\gg=\{(\p(Y)\sub\p(X));(Y\sub X)\in\ovm{\cc}(V_D{}^1),Y\sub V_{D-1}\}$ if $k=0$

\subhead 6.3\endsubhead
If $G$ is of type $A$, then 0.6(b) and 0.7(a) are obvious.

Until the end of 6.9 we assume that $G$ is as in 0.6(iii). In this case, 0.6(b) and
0.7(a) can be verified using the tables in \S5.
We will describe the sets $\ovm{\HH}_c^\gg$ in several examples.

\subhead 6.4\endsubhead
Assume first that $\gg$ is the image in $\un\PP$ of $(P,R)\in\PP$  were $P=G$ hence $R$
is unipotent cuspidal. We will specify $\gg$ by writing the corresponding element of $M(\cg_c)$.
In the following tables
we describe the one element sets $\ovm{\HH}_c^\gg$ for various such $\gg$.

Assume that $|c|=17$ so that $G$ is of type $E_8$. The table is:

$\gg\lra(g_5,\z^j)$:  $(1,H_{5!})$  ($j=1,2,3,4$);

$\gg\lra(g_4,i^j)$:    $(1,H_{4!})$  ($j=1,3$);  

$\gg\lra(g_6,-\th^j)$: $(1,H_{2!3!})$ ($j=1,2$);

$\gg\lra(g_3,\e\th^j)$:     $(1,H_{3!})$ ($j=1,2$);

$\gg\lra(g'_2,\e)$:      $(1,H_{2!2!})$;

$\gg\lra(g_2,-\e)$:      $(1,H_{2!})$;

$\gg\lra(1,\l^4)$:         $(1,1)$.

Assume that $|c|=11$ so that $G$ is of type $F_4$. The table is:

$\gg\lra(g_4,i^j)$:                   $(1,H_{4!})$ ($j=1,3$);

$\gg\lra(g_3,\th^j)$:                   $(1,H_{3!})$  ($j=1,2$);

$\gg\lra(g'_2,\e)$:                      $(1,H_{2!2!})$;

$\gg\lra(g_2,\e)$:                                    $(1,H_{2!})$;

$\gg\lra(1,\l^3)$:                                    $(1,1)$.

Assume that $|c|=4$ so that $G$ is of type $G_2$. The table is:

$\gg\lra(g_3,\th^j)$:   $(1,H_{3!})$ ($j=1,2$);

$\gg\lra(g_2,\e)$:            $(1,H_{2!})$;

$\gg\lra(1,\e)$:             $(1,1)$.

\subhead 6.5\endsubhead
Assume now that $\gg$ is the image in $\un\PP$ of $(P,R)\in\PP$  where
the adjoint group of $\bP$ is of type $E_6$ and $R$ is one of the two unipotent cuspidal
representations of $\bP(F_q)$. In this case $G$ is of type $E_6,E_7$ or $E_8$.

If $|c|=17$ then
$$\ovm{\HH}_c^\gg=\{(1,H_{2!3!}),(H_{2!},H_{2!3!})\}.$$

If $|c|=5$ then
$$\ovm{\HH}_c^\gg=\{(1,H_{3!})\}.$$

\subhead 6.6\endsubhead
Assume now that $\gg$ is the image in $\un\PP$ of $(P,R)\in\PP$  where
the adjoint group of $\bP$ is of type $D_4$ and $R$ is unipotent cuspidal.
In this case $G$ is of type $E_6,E_7$ or $E_8$.

If $|c|=17$ then
$$\ovm{\HH}_c^\gg=
\{(1,H_{2!2!}),(H_{2!},H_{2!2!}),(H_{2!2!},H_8),(H_{3!},H_{2!3!}),(1,H_{2!3!})\}.$$

If $|c|=3$ then
$$\ovm{\HH}_c^\gg=\{(1,H_{2!})\}.$$

\subhead 6.7\endsubhead
Assume now that $\gg$ is the image in $\un\PP$ of $(P,R)\in\PP$  where
the adjoint group of $\bP$ is of type $B_2$ and $R$ is unipotent cuspidal.
In this case $G$ is of type $F_4$.

If $|c|=11$ then
$$\ovm{\HH}_c^\gg=\{(H_{2!},H_{2!2!}),(H_{2!2!},H_8),(1,H_{2!2!})\}.$$

If $|c|=3$ then
$$\ovm{\HH}_c^\gg=\{(1,H_{2!})\}.$$    

\subhead 6.8\endsubhead
Assume now that $\gg$ is the image in $\un\PP$ of $(P,R)\in\PP$  where the adjoint group of
$\bP$ is of type $E_7$ and $R$ is one of the two unipotent cuspidal
representations of $\bP(F_q)$. In this case $G$ is of type $E_7$ or $E_8$.

If $|c|=2$ then $\ovm{\HH}_c^\gg$ is either $\{(1,H_{2!})\}$ or $\{(1,1)\}$.

\subhead 6.9\endsubhead
Assume now that $\gg$ is the image in $\un\PP$ of $(P,R)\in\PP$  where
$\bP$ is a torus and $R$ is the unit representations of $\bP(F_q)$.

If $|c|=17$ then
$$\align&\ovm{\HH}_c^\gg=\ovm{\HH}_c-\{(1,1)\}\\&=
\{(1,H_{2!}),(H_{2!},H_{2!}),(1,H_{2!2!}),(H_{2!},H_{2!2!}),(H_{2!2!},H_{2!2!}),(1,H_{3!}),\\&
(H_{3!},H_{3!}),(H_{2!2!},H_8),(H_8,H_8),(1,H_{2!3!}),(H_{2!},H_{2!3!}),(H_{3!},H_{2!3!}),\\&
(H_{2!3!},H_{2!3!}),(1,H_{4!}),(H_{4!},H_{4!}),(1,H_{5!}),(H_{5!},H_{5!})\}.\endalign$$

If $|c|=11$ then
$$\align&\ovm{\HH}_c^\gg=\ovm{\HH}_c-\{(1,1)\}\\&=
\{(1,H_{2!}),(H_{2!},H_{2!}),(1,H_{2!2!}),(H_{2!},H_{2!2!}),(H_{2!2!},H_{2!2!}),(1,H_{3!}),\\&
(H_{3!},H_{3!}),(H_{2!2!},H_8),(H_8,H_8),(1,H_{4!}),(H_{4!},H_{4!})\}.\endalign$$

If $|c|=4$ then
$$\ovm{\HH}_c^\gg=\ovm{\HH}_c-\{(1,1)\}=
\{(1,H_{2!}),(H_{2!},H_{2!}),(1,H_{3!}),(H_{3!},H_{3!})\}.$$

If $|c|=2$ then
$$\ovm{\HH}_c^\gg=\ovm{\HH}_c-\{(1,1)\}=\{(1,H_{2!}),(H_{2!},H_{2!})\}.$$

If $|c|=5$ then
$$\ovm{\HH}_c^\gg=\ovm{\HH}_c=
\{(1,1),(1,H_{2!}),(H_{2!},H_{2!}),(1,H_{3!}),(H_{3!},H_{3!})\}.$$

If $|c|=3$ then
$$\ovm{\HH}_c^\gg=\ovm{\HH}_c=\{(1,1),(1,H_{2!}),(H_{2!},H_{2!})\}.$$

If $|c|=1$ then
$$\ovm{\HH}_c^\gg=\ovm{\HH}_c=\{(1,1)\}.$$

\head 7. The third basis of $\CC[M(\cg_c)]$\endhead
\subhead 7.1\endsubhead
Until the end of 7.5 we assume that $D$ is even. Let $E\in\cf(V_D)$. Let $\D=D-2\dim(E)$.
For $k\in[0,\D/2]$ let $\fT^k_E$ be as in 2.6 (a subspace of $\fT_E$, see 2.4).
This subspace is in $\pmb\cf(\fT_E)$ hence we can consider its image
$v^k_E\in\fT_E$ under the bijection $\pmb\cf(\fT_E)@>>>\fT_E$ obtained from 3.5(a) by
replacing $V_D$ by $\fT_E$. By \cite{L20a, 1.9} (for $\fT_E$ instead of $V_D$), we have 

$v^k_E=0$ if $k=0$,

$v^k_E=e_{I_{[1,M]}}$ if $k=1$,

$v^k_E=e_{I_{[2,M-1]}}$ if $k=2$,

$v^k_E=e_{I_{[1,2]}}+e_{I_{[M-1,M]}}$ if $k=3$,

$v^k_E=e_{I_{[2,3]}}+e_{I_{[M-2,M-1]}}$ if $k=4$,

$v^k_E=e_{I_{[1,2]}}+e_{I_{[5,M-4]}}+e_{I_{[M-1,M]}}$ if $k=5$,

$v^k_E=e_{I_{[2,3]}}+e_{I_{[6,M-5]}}+e_{I_{[M-2,M-1]}}$ if $k=6$,
\nl 
etc. (Notation of 2.6,)

Assuming that $D\ge2$ and that $i\in[1,D]$, $E'\in\cf(V_{D-2})$ are
such that $E=T_i(E')\op\FF e_i$, we note that (by 2.5(a)) we have for $k\in[0,\D/2]$ 

(a) $T_i(v^k_{E'})=v^k_E$.
\nl
Let $\tu:V_D@>>>\NN$ be as in 3.9; when $D\ge2$ we denote by $\un{\tu}:V_{D-2}@>>>\NN$ the
analogous function. We show:

(b) {\it Let $v\in E$, $k\in[0,\D/2]$,
$z\in\fT^k_E$. If $z=v^k_E$, then $\tu(v+z)=k$. If $z\ne v^k_E$, then $\tu(v+z)<k$.}
\nl
We argue by induction on $\dim(E)$.
If $E=0$ then $v=0$ and the result follows from \cite{L20, 1.15(a)}.
Assume now that $E\ne0$. Then $D\ge2$ and we can find $i\in[1,D]$, $E'\in\cf(V')$
such that $E=T_i(E')\op\FF e_i$. We have $v=T_i(v')+ce_i$ where $v'\in E',c\in\FF$.
We have $z=T_i(z')$ where $z'\in\fT^k_{E'}$; moreover we have $z=v^k_E$ if and
only if $z'=v^k_{E'}$ (see (a)).

We have $\tu(v+z)=\tu(T_i(v'+z')+ce_i)=\un{\tu}(v'+z')$ (we have used
\cite{L20, 1.11(b)}). By the induction hypothesis we have
$\un{\tu}(v'+z')=k$ if $z'=v^k_{E'}$ (that is, if $z=v^k_E$)
and $\un{\tu}(v'+z')<k$ if $z'\ne v^k_{E'}$ (that is, if $z\ne v^k_E$). Now (b) follows.

\subhead 7.2\endsubhead
The following is a reformulation of 3.9(a):

(a) {\it Let $k\in\NN$. Then $v\in V_D$ satisfies $v\in\tu\i(k)$ if and only if
$v=\bar\e_D(\EE)$ where $\EE=E(k)$ with $(E,k)\in\un{\pmb\cf}(V_D)$.} 
\nl
Let $\CC[V_D]$ be the vector space of formal $\CC$-linear combinations of elements of $V_D$.
For $f\in\CC[V_D]$ we write
$f=\sum_{v\in V_D}f(v)v\in\CC[V_D]$ where $f(v)\in\CC$;
let $supp(f)=\{v\in V_D;f(v)\ne0\}$.
For $(E,k)\in\un{\pmb\cf}(V_D)$ we set $f^k_E=\sum_{v\in E+v^k_E}v\in\CC[V_D]$.
For $k\in\NN$ let $\CC[V_D]_k$ (resp. $\CC[V_D]_{\le k}$ be the subspace
$\{f\in\CC[V_D];supp(f)\sub\tu\i(k)\}$ (resp.
$\{f\in\CC[V_D];supp(f)\sub\cup_{k'\in[0,k]}\tu\i(k')\}$) of $\CC[V_D]$. We show:

(b) {\it For any $k\in\NN$, $\{f^k_E;E\in\cf(V_D)\}$ is a $\CC$-basis of $\CC[V_D]_k$.}

(c) {\it $\{f^k_E;(E,k)\in\un{\pmb\cf}(V_D)\}$ is a $\CC$-basis of $\CC[V_D]$.}
\nl
For $\EE\in\pmb\cf(V)$ let $g_\EE=\sum_{v\in\EE}v\in\CC[V_D]$.
By 7.1(b), if $\EE,E,k$ are as in (a), then $g_\EE\in\CC[V_D]_{\le k}$.
By \cite{L20, 1.17(a)}, $\{g_{\EE};\EE\in\pmb\cf(V_D)\}$ is a $\CC$-basis of $\CC[V_D]$.
Hence for any $k\in\NN$, 

(d) $\{g_{\EE};\EE\in\pmb\cf(V_D),\EE=E(k'),(E,k')\in\un{\pmb\cf}(V_D),k'\le k\}$
\nl
is a linearly independent subset of $\CC[V_D]_{\le k}$ of cardinal equal to
$|\cup_{k'\in[0,k]}\tu\i(k')|$ (we use (a)). It follows that

(e) {\it the elements (d) form a $\CC$-basis of $\CC[V_D]_{\le k}$.}
\nl
We show by induction on $k\in\NN$ that

(f) {\it the elements $\{f^{k'}_E;(E,k')\in\un{\pmb\cf}(V_D),k'\in[0,k]\}$ form a $\CC$-basis of
$\CC[V_D]_{\le k}$.}
\nl
Assume first that $k=0$. If $E\in\cf(V_D)$, we have $g_E=f^0_E$ so that (f) follows from (e).
Next we assume that $k\ge1$. Let $(E,k')\in\un{\pmb\cf}(V_D),k'=k$. We have
$g_{E(k)}-f^k_E
=\sum_{v\in E(k)-(E+v^k_E)}v\in\CC[V_D]$; 
moreover, by 7.1(b), $E(k)-(E+v^k_E)$ is contained in $\cup_{k'\le k-1}\tu\i(k')$.
Thus, $g_{E(k)}-f^k_E\in\CC[V_D]_{\le k-1}$. Using the induction hypothesis we see that
$g_{E(k)}-f^k_E$ is a linear combination of elements
$f^{k'}_{E'}$ with $(E',k')\in\un{\pmb\cf}(V_D),k'\le k-1$. 
Thus $g_{E(k)}$ is equal to $f^k_E$ plus a linear combination of elements $f^{k'}_{E'}$ with
$(E',k')\in\un{\pmb\cf}(V_D),k'<k$.
The same is true if here $k$ is replaced by $k_1\in[0,k-1]$ (we use again the
induction hypothesis). It follows that 
the elements $f^{k'}_E$ are related to the elements $g_{E(k')}$ be an upper
triangular matrix with $1$ on diagonal. Hence
(f) follows from (e). By 7.1(b), for any $(E,k)\in\un{\pmb\cf}(V_D)$ we have
$f^k_E\in\CC[V_D]_k$. Hence (b),(c) follow from (f). 

Note that we have now proved 0.6(c) in the case 0.6(i).

\subhead 7.3\endsubhead
Let $\EE,E,k$ be as in 7.2(a). We write $\EE=E(k)$ with $(E,k)\in\un{\pmb\cf}(V_D)$. From 7.1(b)
we see that $\EE\sub\cup_{k'\in[0,k]}\tu\i(k')$ and setting $\EE_*=\{v\in\EE;\tu(v)=k\}$, we have

(a) $\EE_*=E+v^k_E$.

\subhead 7.4\endsubhead
Let $k\in[0,D/2]$ and let $\EE=\la e_{[1,D]}, e_{[2,D-1]},\do,e_{[k,D+1-k]}\ra$
(a primitive subspace of $V_D$). Let $v^k_{\{0\}}=\bar\e_D(\EE)\in V_D$ (a special case of
a notation in 7.1). If $k<D/2$ let $E\in\cf(V_D)$ be the subspace with basis
$e_s,e_{s+2},e_{s+4},\do,e_{D-s}$ where $s=D/2-k$. If $k=D/2$ let $E=0$. Then
$v^k_E\in\fT_E\sub V_D$ is defined as in 7.1. We have the following result.

(a) $v^k_{\{0\}}=v^k_E$.
\nl
Now $v^k_E$ is computed in 7.1 in terms of $I_*$ associated to $E$;
$v^k_{\{0\}}$ is also computed in 7.1, this time in terms of $I_*=(\{1\},\{2\},\do,\{D\})$.
The two computations give the same result; (a) follows.

\subhead 7.5\endsubhead
Assume that $v\in V_D,\tu(v)=k\in[0,D/2]$. We show:

(a) {\it There exists $E'\in\cf(V)$ such that $\dim(E')=D/2-k$ and $v\in E'+v^k_{E'}$.}
\nl
We argue by induction on $D$. Let $\EE,E$ be associated to $v$ as in 7.2(a). By 3.5(b) we
have $v\in\EE$. By the definition of $\EE_*$ (see 3.3) we have $v\in\EE_*$ hence by 7.3(a) we
have $v\in E+v^k_E$.
Assume first that $E=0$. By definition (see 7.1) we have $v=v^k_{\{0\}}$.
Using 7.4(a) we deduce that $v=v^k_{E'}$ where $E'\in\cf(V_D)$ satisfies $\dim(E')=D/2-k$.
Thus (a) holds in this case.

Next we assume that $E\ne0$. Then $D\ge2$, $k\le(D-2)/2$ and we can find
$\EE'\in\pmb\cf(V_{D-2})$ and $i\in[1,D]$ such that $\EE=T_i(\EE')\op\FF e_i$.
Since $v\in\EE$, we have $v=T_i(v')+ce_i$ where $v'\in V_{D-2}$ and $c\in\FF$.
From \cite{L20, 1.11(b)} we have $\tu(v)=\un{\tu}(v')$ (with $\un{\tu}$ as in 7.1).
Thus we have $\un{\tu}(v')=k$. By the induction hypothesis we can find
$E''\in\cf(V_{D-2})$ such that $\dim(E'')=(D-2)/2-k$ and $v'\in E''+v^k_{E''}$.
Let $E'=T_i(E'')+\FF e_i$. We have $E'\in\cf(V_D)$, $\dim(E')=\dim(E'')+1=D/2-k$.
Using 7.1(a) we have

$v=T_i(v')+ce_i\in T_i(E''+v^k_{E''})+\FF e_i\sub E'+v^k_{E'}$.
\nl
Thus $E'$ is as desired. This completes the proof of (a).

\mpb

We see that

(b) {\it $\tu\i(k)$ is a union of affine subspaces of the form $E+v^k_E$ (which are
translates of linear subspaces $E\in\cf(V_D)$ of dimension $D/2-k$).}
\nl
When $k=0$, $\tu\i(k)$ is union of the linear subspaces $E$ in $\cf_*(V_D)$ corresponding to
the various left cell representations. For $k\ge0$, the affine subspaces in (b) can be
regarded as higher left cell representations (generalizing the left cell representations).
When $D=4$ these affine subspaces are as follows:

$(0,1,3,13),(0,2,13,123),(0,1,4,14),(0,3,24,234),(0,2,4,24)$  (if $k=0$),

$(12,124),(124,1234),(1234,134),(134,34)$ (if $k=1$),

$(23)$ (if $k=2$).

(We specify a subset of $V_D$ by a list of its vectors; we write
$i_1i_2\do i_t$ instead of $e_{i_1}+e_{i_2}+\do+e_{i_t}$).

\subhead 7.6\endsubhead
We now assume that $D$ is odd. Most of the results in 7.1-7.5 have
analogues for $V'_D$.
Let $\tu':V'_D@>>>\NN$ be as in 3.10. Then for $k\in\{0,1,3,5,7,\do\}$,
$\tu'{}\i(k)$ is a union of affine subspaces of $V'_D$ which are
translates of linear subspaces $E\in\cf(V'_D)$ of dimension $(D-1)/2-k'$
where $k'=0$ if $k=0$ and $k'=(k+1)/2$ if $k=1,3,5,\do$. When $D=5$ these
affine subspaces are:

$(0,1,3,13),(0,2,13,123),(0,1,4,14),(0,3,24,234),(0,2,4,24)$  (if $k=0$),

$(12,124),(124,1234),(1234,134),(134,34),(23,12)$ (if $k=1$).
\nl
(We specify a subset of $V'_D$ by a list of its vectors; we write
$i_1i_2\do i_t$ instead of $\p(e_{i_1}+e_{i_2}+\do+e_{i_t})$).

\subhead 7.7\endsubhead
Most of the results in 7.1-7.5 have analogues for $M(\cg_c)$ in the setup of 0.1.
Assume for example that $\cg_c=S_5$. The elements $f_{x,\r}$ in 0.6(c) are as follows.

If $(x,\r)\in M(\cg_c)^\gg$ with $\gg$ as in 6.4 then $f_{x,\r}=(x,\r)$.

If $(x,\r)\in M(\cg_c)^\gg$ with $\gg$ as in 6.5 that is,
$(x,\r)\in\{(g_3,\th^j),(g_6,\th^j)\}$ (with $j=1$ or $2$)
\nl
then $f_{x,\r}$ is $(g_3,\th^j),(g_6,\th^j)+(g_3,\th^j)$ respectively.
(The last of these is an analogue of a higher left cell in 7.5.)

If $(x,\r)\in M(\cg_c)^\gg$ with $\gg$ as in 6.6 that is,

$(x,\r)\in\{(g_2,-1),(g_2,-r),(g'_2,r),(g_4,-1)\}$
\nl
then $f_{x,\r}$ is

$(g_2,-1),(g_2,-r)+(g_2,-1),(g'_2,r)+(g_2,-r)+(g_2,-1),(g_4,-1)+(g'_2,r),$

$(g_6,-1)+(g'_2,r)+(g_2,-1)$
\nl
respectively. (The last three of these are analogues of the higher left cells in 7.5.)

If $(x,\r)\in M(\cg_c)^\gg$ with $\gg$ as in 6.9 then $f_{x,\r}=g_{x,\r}$ (notation of 0.6).

\enddocument